\documentstyle[11pt]{article}
\input amssym.def
\input amssym

\newtheorem{theorem}{Theorem}[section]
\newtheorem{lemma}[theorem]{Lemma}
\newtheorem{proposition}[theorem]{Proposition}

\newtheorem{conjecture}[theorem]{Conjecture}
\newtheorem{definition}[theorem]{Definition}

\begin{document}
\begin{titlepage}
\title{ \bf Deninger's conjecture on $L$-functions of elliptic curves
  at $s=3$.}

\author{Alexander Goncharov}
\end{titlepage}
\maketitle
\date{}

\tableofcontents

\section  {Introduction}
\vskip 3mm \noindent
{\bf 1.} In this paper I  compute explicitly the
regulator map on $K_4(X)$
for an arbitrary curve $X$ over a number field. Using this
and Beilinson's
theorem about regulators for modular curves ([B2]) I prove a
formula expressing the value of the $L$-function $L(E,s)$ of a
modular elliptic curve $E$ over $\Bbb Q$ at $s=3$ by the double
Eisenstein-Kronecker series. It was
conjectured by C. Deninger [D1].

{\bf 2. Generalized Eisenstein-Kronecker series}.
Let $E$ be an elliptic curve  and $\Gamma :=
H_1(E(\Bbb C),\Bbb Z)$. Choose a
 holomorphic
1-form $\omega$. It defines an embedding
$\Gamma \hookrightarrow \Bbb C$ together with an isomorphism
$E(\Bbb C) = \Bbb C/\Gamma  =  \Gamma \otimes\Bbb R/\Gamma$.
The Poincare duality provides a nondegenerate pairing $\Gamma\times
 \Gamma  \longrightarrow \Bbb Z(1)$.
 Let
$$
(\cdot ,\cdot ): E(\Bbb C)\times \Gamma \longrightarrow
\Bbb R(1)/\Bbb Z(1) = U(1) \subset \Bbb C^{\ast}
$$
be the corresponding pairing.
If $\Gamma = \Bbb Zu + \Bbb Zv \subset \Bbb C$  with
$Im(u/v) >0$ then $(z,\gamma) = \exp A(\Gamma)^{-1}(z\bar \gamma
-\bar z \gamma)$ where
 $A(\Gamma) = \frac{1}{2\pi i}(\bar u v-u \bar v)$.

Let $x,y,z \in E(\Bbb C)$ and $n \geq 3$.
The function
\begin{equation}
K_n(x,y,z) := \sum'_{\gamma_1+...+\gamma_n=0}
\frac{(x,\gamma_1)(y,\gamma_2 + ...+\gamma_{n-1})(z,\gamma_n)(\bar\gamma_n -
\bar\gamma_{n-1})}{\vert\gamma_1\vert^2 \vert\gamma_2\vert^2 ...
\vert\gamma_n\vert^2}
\end{equation}
will be called generalized Eisenstein-Kronecker series. It
is invariant under the shift $(x,y,z) \rightarrow (x+t,y+t,z+t)$ and so
lives actually  on  $E(\Bbb C) \times E(\Bbb C)$.

To formulate the results I have first to recall the definition of

{\bf 3. The group $B_2(F)$}.
Let $F$ be a field and $\Bbb Z[F^*]$ be
the free abelian group generated by symbols $\{x\}$ where $x \in F^*$.
Let $R_2(F)$ be the subgroup of $\Bbb Z[F^*]$  generated by
the elements
\begin{equation} \label{1}
\sum_{i=1}^5(-1)^i \{r
(x_{1},\ldots ,\hat x_{i},\ldots ,x_{5})\}
\end{equation}
where $x_{1},\ldots ,x_{5}$
run through all 5-tuples of distinct $F$-points  of $P^{1}$.
By definition
$
B_2(F) := \Bbb Z[F^*]/R_2(F)
$.
Symbol $\{x\}_2$ denotes the projection of  $\{x\}$ to
$B_2(F)$.

One can show that the formulas
%\begin{equation}   \label{0q}
 $\delta:  \{x\}
\longmapsto
(1-x)\wedge x; \qquad  \{1\} \longmapsto 0$
%\end{equation}
provide us with a homomorphism of groups
$\delta: B_2(F) \longrightarrow \wedge^2F^*$. (In other words
$\delta(R_2(F)) = 0$).
This is  one of the most important properties of the group $B_2(F)$.

{\bf 4. Special values of L-functions}. Let $v_x : \Bbb Q(E)^{\ast}
\longrightarrow \Bbb Z$ be the valuation
defined by a point $x \in E(\bar \Bbb Q)$. Denote by $f_E$  the conductor of
$E$.
Let $\omega \in
H^0(E,\Omega^1_{E/\Bbb R})$,
Denote by $
\Omega = \int_{E(\Bbb R)}\omega$
the real period of $E$.

\vskip 3mm \noindent
\begin {theorem}  \label {0.1}
 Let $E$ be a modular elliptic curve over $\Bbb Q$.
Then there exist rational functions $f_i,g_i \in \Bbb Q(E)^{\ast}$
satisfying the conditions:
\begin{equation}   \label{c3qm}
\sum_i(1-f_i)\wedge f_i \wedge g_i = 0 \quad \mbox{in} \quad
\Lambda^3 \Bbb Q(E)^{\ast}
\end{equation}
\begin{equation}   \label{o1}
 \sum_i
v_x(g_i)\{f_i(x)\}_2 = 0 \quad \mbox{in} \quad
B_2(\bar \Bbb Q) \quad \mbox{for any} \quad x \in E(\bar \Bbb Q)
\end{equation}
such that
\begin{equation}    \label{c4q}
L(E,3) = q \Bigl(\frac{2 \pi A(\Gamma)}{f_E} \Bigr)^2 \Omega \cdot
\sum_i\sum'_{\gamma_1+\gamma_2 +\gamma_3=0}
\frac{(x_i,\gamma_1)(y_i,\gamma_2 )(z_i,\gamma_3)(\bar\gamma_3 -
\bar\gamma_2)}{\vert\gamma_1\vert^2 \vert\gamma_2\vert^2
\vert\gamma_3\vert^2 }
\end{equation}
 where $q$ is a non-zero rational number and $x_i,y_i,z_i$ are the
divisors of the functions $g_i, f_i, 1-f_i$ respectively.
\end {theorem}
It is interesting that the right-hand side of (\ref{c4q}) depends
only on the divisors of the functions $g_i, f_i, 1-f_i$.

A similar formula expressing $ L(E,2)$  for a modular elliptic curve over $\Bbb
Q$
by the classical Eisenstein-Kronecker series
$$
L(E,2) =  q \Bigl(\frac{2 \pi A(\Gamma)}{f_E} \Bigr) \Omega \cdot\sum_i
\sum'_{\gamma}
\frac{(x_i-y_i,\gamma)\bar\gamma}
{\vert\gamma\vert^4}
$$
was known thanks to Bloch and Beilinson [Bl1], [B1].

A  formula (\ref{c4q})  for an arbitrary
elliptic curve $E$ over $\Bbb Q$ in a slightly different form  was
conjectured
by C.Deninger ([D1]), who used Massey products in Deligne
cohomology to guess a formula for $L(E,3)$.  I do not use
Massey products  in the formulation or proof of the theorem.

One can define ([G2]) for an arbitrary
field $F$  an  abelian group
$$
{\cal B}_{n}(F):=\Bbb Z[P^{1}_{F}]/{\cal R}_{n}(F)
$$
together with a homomorphism
$$
{\cal B}_n(F) \stackrel{\delta}{\longrightarrow} B_{n-1}(F) \otimes F^{\ast}
\qquad \{x\}_n
\longmapsto
\{x\}_{n-1}\otimes x
$$
I will recall the  definition of ${\cal R}_{n}(F)$ in chapter 4 below.
Roughly speaking it is the ``connected component of zero'' of  ${\rm
Ker} \delta$.
One can show that ${\cal R}_{n}(\Bbb C)$  is the subgroup of all functional
equations for the
classical $n$-logarithm, see [G2].

\vskip 3mm \noindent
\begin {conjecture}
\label {0.3}
  Let $E$ be an elliptic curve over $\Bbb Q$.
Then
there exist rational functions $f_i,g_i \in \Bbb Q(E)^{\ast}$
satisfying the condition ($n > 3$)
\begin{equation} \label{c1}
\sum_i\{f_i\}_{n-2}\otimes f_i \wedge g_i \in
 {\cal B}_{n-2}(\Bbb Q(E)) \otimes  \Lambda^2  \Bbb Q(E)^{\ast}
\end{equation}
\begin{equation}   \label{c2}
 \sum_i
v_x(g_i)\{f_i(x)\}_{n-1} = 0 \quad \mbox{in} \quad
{\cal B}_{n-1}({\bar \Bbb Q}) \quad \mbox{for any} \quad x \in
E(\bar \Bbb Q)
\end{equation}
such that
\begin{equation}  \label{c2222}
q \cdot L(E,n) =  \Bigl(\frac{2 \pi A(\Gamma)}{f_E} \Bigr)^{n-1}
\Omega  \cdot \sum_i K_n(x_i,y_i,z_i)
\end{equation}
where $q$ is a non-zero rational number and $x_i,y_i,z_i$ are the divisors
of the functions $g_i, f_i, 1-f_i$.
\end {conjecture}

\begin{conjecture}
 For any $f_i,g_i \in \Bbb Q(E)^{\ast}$ satisfying the conditions
 of theorem (\ref{0.1}) (resp. conjecture (\ref{0.3}))  one has (\ref{c2222})
 with
 $q \in \Bbb Q$.
\end{conjecture}

Beilinson's  conjecture on $L$-functions  permits to
formulate  a similar conjecture for an elliptic curve over any
number field $F$ in which we have in the right hand side a
determinant  whose entries are the functions $K_n(x,y,z)$.

{\bf Remark}. If $n>2$ then for a regular proper model $E_{\Bbb
Z}$ of $E$ over $Spec(\Bbb Z)$ one has $gr^{\gamma}_nK_{2n-2}(E_{\Bbb Z}) =
gr^{\gamma}_nK_{2n-2}(E)$, so, unlike to the $n=2$ case,  we
don't have to  worry
about the ``integrality condition''.

{\bf 5. Explicit formulas for the regulators for curves}.
 Let $X$ be a curve  over $\Bbb R$
and $n>1$.
Then the  real Deligne cohomology $H_{{\cal D}}^2(X/\Bbb R,\Bbb R(n))$
  equals
$H^1(X/\Bbb R,\Bbb R(n-1))$. Further,
cup product with $\omega \in \Omega^1(\bar X)$ provides an isomorphism
of vector spaces over $\Bbb R$:
$$
H^1(X/\Bbb R,\Bbb R(n-1)) \longrightarrow H^0(\bar X, \Omega^1)^{\vee}
$$
So we will present elements of $H_{{\cal D}}^2(X/\Bbb R,\Bbb R(n))$ as
functionals on $H^0(\bar X, \Omega^1)^{\vee}$.

In chapter 3  we prove  the following explicit formulas for the regulators
$$
r_{{\cal D}}(3) : K_4(X) \longrightarrow H_{{\cal D}}^2(X/\Bbb R,\Bbb R(3))
$$
which generalize the famous symbol  on $K_2(X)$ of Beilinson and
Deligne ([B3], [Del]). (For simplicity we formulate results for
curves over $\Bbb Q$).

I will use the notation
\begin {equation} \label{al11}
\alpha(f,g) : = \log \vert  f\vert  d \log \vert g\vert  -
\log \vert  g\vert d \log \vert f\vert
\end {equation}
\begin {theorem}
\label {ohti}
Let $X$ be a regular  curve
over $\Bbb Q$. Then for each element $\gamma_4 \in K_4(X)$ there are
rational functions
$f_i,g_i \in \Bbb Q(X)$ satisfying the conditions (\ref{c3qm}) and (\ref{o1})
such that  for any $\omega \in \Omega^1(\bar X)$ one has
$$
\int_{X(\Bbb C)} r_{{\cal D}}(3)(\gamma_4)\wedge \omega =
\int_{X(\Bbb C)}
 \log|g_i|\alpha(1-f_i,f_i)\wedge \omega
$$
\end {theorem}

The proof of theorem (\ref{ohti}) is based on the results of [G2],
[G3].

The regulator map on a certain subgroup of $K_4^{(3)}$ of curves
over number fields was also computed by R. de Jeu [J].

In general the Beilinson regulator is a map
$$
r_{{\cal D}}(n+1) : K_{2n}(X) \longrightarrow H_{{\cal D}}^2(X/\Bbb R,\Bbb
R(n+1))
$$
We expect the following to be true.
\begin {conjecture} \label{vbnm}
Let $X$ be a nonsingular  curve
over $\Bbb Q$. Then for each element $\gamma_{2n} \in K_{2n}(X)$
there are rational functions
$f_i,g_i \in \Bbb Q(X)$ satisfying the conditions (\ref{c1}) and (\ref{c2})
such that  for any $\omega \in \Omega^1(\bar X)$
$$
\int_{X(\Bbb C)} r_{{\cal D}}(n+1)(\gamma_{2n})\wedge \omega =
c_{n+1}\cdot \sum_i\int_{X(\Bbb C)}
 \log|g_i| \log^{n-2}|f_i|\alpha(1-f_i,f_i)\wedge \omega
$$
where $c_{n+1} \in \Bbb Q^{\ast}$ is a certain explicitly computable constant.
\end {conjecture}
Moreover, one can prove that condition (\ref{c3qm}),  or respectively
(\ref{c1}) if $n>3$,  implies that the
right-hand side of these formulas depends
only on the divisors of the functions $f_i, g_i, 1-f_i$ .
When $X$ is an elliptic curve   this  together with
 Fourier transform  and Beilinson's theorem ([B2]) lead
to formulas for $L(E,n)$  from theorem 1.1   and
 conjectures 1.2 - 1.3.

The proof of this conjecture for  $K_6(X)$  will be published in
[G4].

In chapter 4 we will see that  conjecture (\ref{vbnm})  follows from the main
conjecture in
[G2]  which tells  us that the  complexes $\Gamma(F,n)$ constructed there
catch all of the
rational algebraic $K$-theory of an {\it arbitrary} field $F$.

The crucial role in the proof of these results is played by the classical
 $n$-logarithms $Li_n(z) =
\int_0^zLi_{n-1}(t) d\log t$. The single-valued version of the
$n$-logarithm is the following function ([Z1]):
\begin{eqnarray*}
{\cal L}_{n}(z) &:=& \begin{array}{ll}
{\rm Re} & (n:\ {\rm odd}) \\
{\rm Im} & (n: \ {\rm even}) \end{array}
\left( \sum^{n}_{k=0} \beta_k
\log^{k}\vert z\vert \cdot Li_{n-k}(z)\right)\; , \quad n\geq
2 \\
\end{eqnarray*}
Here   $\beta_k
=B_{k}\cdot 2^k/k!$  and  $B_k$ are Bernoulli numbers:
$\sum_{k=0}^{\infty}\beta_k x^k = \frac{2x}{e^{2x} -1}$.

One can show (see proposition (\ref{ma1})) that for any functions $f_i,g_i$
satisfying (\ref{c2}) one can write the regulator integrals
$$
\int_{X(\Bbb C)}
 \log|g_i| \log^{n-2}|f_i|\alpha(1-f_i,f_i)\wedge \omega =
b_{n+1}\cdot\sum_i\int_{X(\Bbb C)} {\cal L}_{n}(f_i)d\log|g_i|\wedge \omega
$$
where  $b_{n+1}$ are certain explicitly computable non zero rational
constant.

{\bf 6. The structure of the paper}.
Let ${\cal O}$ be a local ring with infinite residue field.
In chapter 2 we will construct homomorphisms
\begin{equation}  \label{c+2}
K^{[i]}_{6-i}({\cal O})_{\Bbb Q} \longrightarrow H^i( \Gamma ({\cal O},3))
\end{equation}
where $K^{[i]}_{n}({\cal O})$ are the graded quotients of the rank
filtration on Quillen's K-groups of the
ring ${\cal O}$. Hypothetically modulo torsion it is
opposit to the Adams
filtration.

Now let $X$ be a curve over a number field $F$. Then we define a
complex $ \Gamma(X,3)$ and  homomorphisms
\begin{equation}  \label{c+3}
K^{[i]}_{6-i}(X)_{\Bbb Q} \longrightarrow H^i( \Gamma (X,3))
\end{equation}
We impose  the condition that $F$ is a number field only because of one
argument
``ad hoc'' in proof which is based on the Borel theorem.
It would be interesting to prove this result for an arbitrary field $F$.
This goal is almost  (but not completely) achieved in the sections 4-7 of
chapter 2.
The reader who is
interested only in the proof of the formula
for $L(E,3)$ can skip it.

In chapter 3 we will prove that the composition of this map
with the natural map from $\Gamma (X,3)$ to Deligne cohomology
coincides with Beilinson's regulator.

In the end of chapter 3 and in chapter 4 we compute the regulator integrals for
curves. In particular we show that conjecture 1.5 essentially
follows from the main conjecture of [G1] on the structure of motivic
complexes.
Then we apply these results to elliptic curves and get
the generalized Eisenstein-Kronecker series. Therefore we finish the proof of
theorem 1.1 and deduce conjecture 1.2 from conjecture 1.5.

{\bf Acknowledgement}. This work was partially supported by the
NSF Grant DMS-9500010. The final version of this paper was completed
during my stay in
MPI (Bonn) whose hospitality and support are gratefully acknowledged.

\section  {  Motivic complexes for a curve and algebraic  $K$-theory  }

{\bf 1 The weight 3  motivic complex }. We will call by this name
the complex $\Gamma(X,3)$ introduced in [G1-G2] for an arbitrary
regular scheme $X$. If $X = Spec(F)$ where $F$ is an {\it
  arbitrary} field, it looks as follows:
\begin{equation}
B_3(F) \longrightarrow B_2(F)\otimes
F^* \longrightarrow \Lambda^3F^*
\end{equation}
Here $B_3(F) := \Bbb
Z[F^*]/R_3(F)$ where the subgroup $R_3(F)$ will be defined below,
after theorem (\ref{z2}).
The group $B_3(F)$ is placed in degree 1. The differential has degree
$+1$ and is defined as follows:
$$
\{x\}_3 \longmapsto \{x\}_2\otimes
x; \quad \{x\}_2 \otimes y \longmapsto (1-x) \wedge x \wedge y
$$

In [G1] we have  constructed
 homomorphisms
of groups
$$
c_{2,3}: K_{4}(F)\otimes \Bbb Q \longrightarrow
H^2\Gamma(Spec(F);3) \otimes \Bbb Q
$$
We will recall the definition of  this homomorphism below.

The goal of this chapter is to get a  similar homomorphism for
curves over number
fields.

{\bf 2. The complex $\Gamma(X;3)$ for a curve $X$ over a field $F$.}
Let $K$ be an arbitrary  field with
 discrete valuation $v$ and  residue class $
k_v$.  The group of units $U$ has a natural
homomorphism $U\longrightarrow k_v^{\ast}\; , \; u\mapsto
\bar u$.  An element $\pi \in K^{\ast}$ is prime if
ord$_{v}\pi = 1$.

Let us define the residue homomorphism
\begin{equation}
\partial_{v} :\Gamma (K,3)\longrightarrow \Gamma( k_v,
2)[-1]
\end{equation}

There is a homomorphism $\theta_n:\wedge^{n}K^{\ast}
\longrightarrow\wedge^{n-1} k_v^{\ast}$ uniquely defined
by the  properties $(u_{i}\in U)$:

 $\theta_n (\pi\wedge u_{1}\wedge \cdots\wedge u_{n-1})
= \bar u_{1}\wedge\cdots \wedge \bar u_{n-1}$ and  $\theta_n
(u_{1}\wedge \cdots \wedge u_{n}) = 0$.

It clearly does not depend on the choice of $\pi$.
Let us define a homomorphism $s_{v}:\Bbb
Z[K^*]\longrightarrow \Bbb Z[ k_v^*]$ as
follows
$$
s_{v}\{ x\} = \left\{
\begin{array}{ll}
\{ \bar x\} & \mbox{ if $x$ is a unit}\\
0 & {\rm otherwise} \end{array}
\right \}
$$

Then it induces a homomorphism
$s_{v} :  B_{2}(K)\longrightarrow B_{2}(k_v)$ (see s.\ 9
\S 1 of [G1]).
We get a homomorphism
$$
s_{v}\otimes \theta_1 : B_{2}(K)\otimes
K^{\ast}\longrightarrow B_{2} (k_v)
$$
Let us consider the following map $\partial_{v}$  of complexes:
\begin{equation}
\begin{array}{ccccc}
B_3(K)&\stackrel{\delta}{\longrightarrow}&B_2(K)\otimes
K^{\ast}&
\stackrel{\delta}{\longrightarrow}&\Lambda^3K^{\ast}\\
  &&\downarrow s_{v}\otimes \theta_1 &&\downarrow  \theta_3 \\
&&B_2(k_v)&\stackrel{\delta}{\longrightarrow}&\Lambda^2k_v^{\ast}
\end{array}
\end{equation}

%\begin{lemma} \label{2.100}
%\end{lemma}

The maps $\partial_{v}$ define a homomorphism of complexes, see s.\ 14 of \S
1 in [G1].

Let $k_x$ be the residue field
at the point $x \in X$.

By definition $\Gamma(X;3)$ is
the total complex associated with the   bicomplex
$$
\begin{array}{ccccc}
B_3(F(X))&\stackrel{\delta}{\longrightarrow}&B_2(F(X))\otimes
F(X)^{\ast}&
\stackrel{\delta}{\longrightarrow}&\Lambda^3F(X)^{\ast}\\
  &&\downarrow \partial&&\downarrow \partial\\
&&\coprod_{x\in X_1}B_2(k_x)&\stackrel{\delta}{\longrightarrow}&\coprod_{x\in
X_1}\Lambda^2k_x^{\ast}
\end{array}
$$

Here $\partial = \coprod_{x\in X_1}\partial_x$, where $\partial_x$ is
the residue homomorphism related to the valuation on $F(X)$
 corresponding to the
point $x$; the very left group  placed in degree 1 and the
differentials
have degree $+1$.

{\bf 3. The key result}. The main point is to
  show that homomorphism $c_{i,3}$
   carries the
  residue map in Quillen K-theory to the one on
  $\Gamma$-complexes.

One has the exact localization sequence
\begin{equation}
\longrightarrow K_{n}(X)\longrightarrow K_{n}(F(X)) \stackrel{\tilde
\delta}{\longrightarrow}
\coprod_{x\in X_1} K_{n-1}(k_x) \longrightarrow
\end{equation}
where $X_1$ is the set of all codimension one points of a scheme $X$
and $\tilde \delta$ is the residue homomorphism in the Quillen
$K$-theory.

So keeping in mind the localization sequence  we see that in
order to construct a homomorphism of groups
$$
K_{4}(X)\otimes \Bbb Q \longrightarrow H^2(\Gamma(X;3) \otimes \Bbb Q)
$$
the only thing we have to prove is the following

\begin{theorem} \label{mth}
Let $F$ be a number field. Then the  diagram
$$
\begin{array}{ccc}
K_{4}(F(X))&\stackrel{c_{2,3}}{\longrightarrow}&H^2\Gamma(Spec(F(X));3)_{\Bbb
Q} \\
&&\\
\downarrow \tilde \delta&&\downarrow  \delta\\
&&\\
K_{3}(k_x)&\stackrel{ c_{1,2}}{\longrightarrow}&
H^{1}\Gamma(Spec(k_x);2)_{\Bbb Q}
\end{array}
$$
is commutative.
\end{theorem}

{\bf Remark}. It is only important for us that $c_{2,3} \circ \delta = q \tilde
\delta c_{1,2}$
for a certain nonzero constant $q \in \Bbb Q^*$.

This theorem will be proved in the section 7 of chapter 3.

It would be interesting to prove this theorem for an arbitrary field $F$. Here
is
a possible strategy.
Let $K_x$ be the completion of the field $F(X)$ at the point $x$.
Denote by ${\cal O}_x$ the ring of integers in $K_x$. Then $k_x$ is the residue
field.
 Let $i: K_4(F(X)) \longrightarrow K_4(K_x)$ be the natural map.
One has the folowing  diagram:
$$
\begin{array}{ccc}
K_{4}(F(X))&\stackrel{c_{2,3}}{\longrightarrow}&H^2\Gamma(Spec(F(X));3)_{\Bbb
Q} \\
&&\\
\downarrow i&&\downarrow  i'\\
K_{4}(K_x)&\stackrel{c_{2,3}}{\longrightarrow}&H^2\Gamma(Spec(K_x);3)_{\Bbb
Q} \\
&&\\
\downarrow \tilde \delta&&\downarrow  \delta\\
&&\\
K_{3}(k_x)&\stackrel{ c_{1,2}}{\longrightarrow}&
H^{1}\Gamma(Spec(k_x);2)_{\Bbb Q}
\end{array}
$$
The group
$K_{4}(K_x)$ is generated by $K_{4}({\cal O}_x)$ and
  $K_{3}(k_x) \cdot K_x^{\ast}$ where $\cdot $ is the product in
Quillen's K-theory.
So to prove the theorem for an arbitrary field $F$ we have to show that

a) {\it The composition
$$
K_{4}({\cal O}_x) \longrightarrow
K_{4}(K_x)\stackrel{c_{2,3}}{\longrightarrow} H^2\Gamma
(Spec(K_x);3)_{\Bbb Q} \stackrel{ \delta}{\longrightarrow}
H^{1}\Gamma(Spec(k_x);2)_{\Bbb Q}
$$
 is zero.

b) The statement of the theorem is true for the subgroup}
$K_{3}(k_x) \cdot K_x^{\ast} \cap i(K_4(F(X)))$.

The rest of this chapter is devoted to the proof of statement a).

{\bf 4. Proof of  a): the beginning.} In this section we suppose that ${\cal
O}$ is a local ring with an infinite residue field $k$.  Let ${\cal
O}^{\infty}$ be the free
${\cal O}$-module with the basis $e_1,...,e_n,...$ and ${\cal O}^{n}$
be the submodule with the basis $e_1,...,e_n$. A vector
$v \in {\cal O}^{n}$ will be identified with the corresponding
column of height $n$. By definition a set of vectors
$v_1,...,v_m \in {\cal O}^{n}$
is  jointly unimodular  if  the matrix $(v_1,...,v_m)$ is left invertible
in $M_{nm}({\cal O})$. Any projective module over ${\cal O}$
is free, so
one can show that any jointly unimodular set of vectors can be
completed to a basis of ${\cal O}^{n}$.

Let $V$ be a free ${\cal O}$-module of rank $n$ and $v_1,...,v_m \in V$.
We will say that the vectors $v_i$ are in general position if
any $min(n,m)$
of them are jointly unimodular.
This notion is independent of choice of a basis in $V$.

Let ${\tilde C}_{k}({\cal O}^{n})$ be the free abelian group generated by
$k+1$-tuples
 of vectors in generic position in ${\cal O}^{n}$. They form a complex
${\tilde C}_{\ast}({\cal O}^{n})$ with the differential $d$ given by the
usual formula
$$
d:  {\tilde C}_{m}({\cal O}^{n})
\to   {\tilde C}_{m-1}({\cal O}^{n});\quad d:(v_{1},\ldots ,v_{m+1})\mapsto
\sum^{m+1}_{i=1}(-1)^{i-1}(v_{1},\ldots ,\hat
v_{i},\ldots,v_{m+1})\; .
$$
This complex is acyclic in degrees bigger then $0$ and so is a
resolution of the trivial $GL_n({\cal O})$-module $\Bbb Z$. The
group $GL_n({\cal O})$ acts on ${\tilde C}_{k}({\cal O}^{n})$.

Configurations of $m$
vectors in ${\cal O}^{n}$ are $m$-tuples of vectors considered modulo
$GL({\cal O}^{n})$-equivalence.

We get a complex $C_{\ast}({\cal O}^{n})$.
One has canonical homomorphism
$$
H_i(GL_n({\cal O}) \longrightarrow H_i(C_{\ast}({\cal O}^{n}))
$$
Set
$$
p_2: \Bbb Z[F^* ] \longrightarrow B_2(F) \qquad p_3:
\Bbb Z[F^* ] \longrightarrow B_3(F)
$$
$$
B_2({\cal O}):= p_2(\Bbb Z[{\cal
  O}^* ])\qquad B_2({\cal O}):= p_3(\Bbb Z[{\cal
  O}^* ])
$$
Then one has complexes
$$
B_2({\cal O}) \longrightarrow \Lambda^2{\cal O}^*
$$
$$
B_3({\cal O}) \longrightarrow  B_2({\cal O})\otimes {\cal
  O}^*  \longrightarrow
\Lambda^3{\cal O}^*
$$
which are subcomplexes of $\Gamma(Spec (K),2)$ and $\Gamma(Spec (K),3)$.
Let us construct a homomorphism of complexes
\begin{center}
\begin{picture}(100,70)(20,0)
\put(-15,50){$C_{5}({\cal O}^{3})$}
\put(0,45){\vector(0,-1){40}}
\put(+5,25){$f_{6}(3)$}
\put(-28,-5){$B_{3}({\cal O})$}
\put(70,50){$C_{4}({\cal O}^{3})$}
\put(85,45){\vector(0,-1){40}}
\put(90,25){$f_{5}(3)$}
\put(-75,54){\vector(1,0){40}}
\put(30,0){\vector(1,0){30}}
\put(60,-5){$ B_{2}({\cal O})\otimes {\cal O}^{\ast}$}
\put(130,0){\vector(1,0){20}}
\put(35,54){\vector(1,0){35}}
\put(110,55){\vector(1,0){40}}
\put(155,50){$C_{3}({\cal O}^{3})$}
\put(155,-5){$\Lambda^{3}{\cal O}^{\ast}$}
\put(170,45){\vector(0,-1){40}}
\put(175,25){$f_{4}(3)$}
\put(260,25){}
\end{picture}
\end{center}
We will use the following notation. For any $n$ vectors $v_1,...,v_n$ in
${\cal O}^{n}$ set $\Delta(v_1,...,v_n):= det (v_1,...,v_n)$
where $det (v_1,...,v_n)$ is the $n\times n$ matrix formed by the columns of
coordinates of vectors $v_i$ in the canonical
basis $e_1 = (1,0,...,0), ..., e_n = (0,...,0,1)$.
Notice that  vectors $v_1,...,v_n$ are in generic position (=
jointly unimodular) if and only if $\Delta(v_1,...,v_n) \in {\cal O}^*$

Let ${\rm Alt}_n\ f(v_{1},\ldots ,v_{n}):= \sum_{\sigma \in S_{n}}(-
1)^{\vert \sigma\vert} f(v_{\sigma(1)},\ldots
,f_{\sigma(n)}$.
Set
$$
f_{4}(3): (v_{1},\ldots ,v_{4})\mapsto {\rm Alt}_4\
\Delta(v_{1},v_{2},v_{3})\wedge \Delta (v_{1},v_{2},v_{4})
\wedge \Delta (v_{1},v_{3},v_{4}) % \; .
$$
$$
f_{5}(3)(v_{1},\ldots ,v_{5}):= \frac{1}{2} {\rm
  Alt}_5(\{r(v_{1}\vert  v_{2},\ldots ,  v_{5})\}_{_{\scriptstyle
2}}\otimes \Delta (v_{1},v_{2} , v_{3}))\; .
$$
Here $(v_{1}\vert  v_{2},\ldots ,  v_{5})$ is the configuration
of four vectors in $V/<v_1>$ obtained by the projection of
vectors $v_2,...,v_5$. We take then
the cross-ratio of the corresponding points on the projective
line.
Now put
\begin{equation} \label{z3}
f_{6}(3):(v_{1},\ldots ,v_{6})\mapsto \frac{1}{15}{\rm Alt}_6\left\{
\frac{\Delta (v_{1},v_{2},v_{4})\Delta (v_{2},v_{3},v_{5})
\Delta(v_{3},v_{1},v_{6})}
{\Delta(v_{1},v_{2},v_{5})\Delta (v_{2},v_{3},v_{6})\Delta
(v_{3},v_{1},v_{4})}\right\}
\end{equation}

\begin{definition} The subgroup $R_3(F) \subset \Bbb Z[F^*]$ is
generated by the elements $\sum_{i=1}^7
(-1)^i f_6(3)(v_1,...,\hat v_i,...,v_7)$.
 \end{definition}

\begin{theorem} \label{z2}
a) $f_{4}(3)$  and $f_{5}(3)$ do not depend on the choice
of $\omega$.

b) The homomorphisms $f_*(3)$ provide a morphism of
complexes.
\end{theorem}

{\bf Proof}. See the appendix.

{\bf 5.  What remains to be done}.
  Just by the construction we have a commutative diagram where
the vertical arrows are
the natural inclusions:
$$
\begin{array}{ccccc}
C_*({\cal O}^3)&\stackrel{f_*(3)}{\longrightarrow}
&\Gamma({\cal O};3)&\stackrel{\delta}{\longrightarrow}&0\\
&&&&\\
\downarrow &&\downarrow &&\downarrow\\
&&&&\\
C_*(K^3)&\stackrel{f_*(3)}{\longrightarrow}&\Gamma(K;3)&\stackrel{\delta}{\longrightarrow}&\Gamma(k;2)\\
\end{array}
$$
Thus we have constructed a homomorphism
\begin{equation} \label{z4}
H_{4}(GL_3({\cal O})) \longrightarrow H^2\Gamma({\cal O};3)
\end{equation}
such that the composition
$$
H_{4}(GL_3({\cal O})) \longrightarrow
H_{4}(GL_3(K))\stackrel{c_{2,3}}{\longrightarrow} H^2\Gamma
(K;3)_{\Bbb Q} \stackrel{ \delta}{\longrightarrow}
H^{1}\Gamma(k;2)_{\Bbb Q}
$$
is equal to zero.

To complete the part a) of our program we have to do the
stabilisation, i.e. for any $n >3$ to extend the homomorphism (\ref{z4}) to a
  homomorphism
$$
H_{4}(GL_n({\cal O})) \longrightarrow H^2\Gamma({\cal O};3)
$$
which fits into a commutative diagram
$$
\begin{array}{ccc}
H_{4}(GL_n({\cal O})) &\longrightarrow &H^2\Gamma({\cal O};3)\\
&&\\
\downarrow &&\downarrow\\
&&\\
H_{4}(GL_n(K)) & \longrightarrow &  H^2\Gamma(K;3)
\end{array}
$$
This will be done in the next three sections.

{\bf 6.  The bi-Grassmannian complex  over a field ([G2])} The
bicomplex
\begin{equation}
\begin{array}{cccccccc}
& && & & \downarrow\qquad  & & \downarrow \quad
\\
 & && & \stackrel{d}{\longrightarrow}
& C_{7}(5) & \stackrel{d}{\longrightarrow} & C_{6}(5)
\\
&&& \downarrow d' & & \downarrow d' & & \downarrow d' \\
&&\stackrel{d}{\longrightarrow}  & C_{7}(4) & \stackrel{d}{\longrightarrow}
 &C_{6}(4) & \stackrel{d}{\longrightarrow} & C_{5}(4)
\\
&\downarrow d'&& \downarrow d'\; \; & & \downarrow d' & & \downarrow d' \\
\stackrel{d}{\longrightarrow}&C_{7}(3)&\stackrel{d}{\longrightarrow} & C_{6}(3)
& \stackrel{d}{\longrightarrow}
& C_{5}(3) & \stackrel{d}{\longrightarrow}& C_{4}(3)
\end{array}
\end{equation}
where
$$
d':(l_{1},\ldots ,l_{m})\mapsto {\displaystyle \sum^{m}_{i=1}}
(-1)^{i-1}(l_{i}\vert l_{1},\ldots ,\hat{l_{i}},\ldots ,
l_{m})
$$
is  the weight three Grassmannian bicomplex.
Denote by $(BC_{\ast}(3),\partial)$ the corresponding total
complex. We place $C_{4}(3)$ in degree 3 and $\partial$ has degree $-1$.
We define a map of complexes $\psi_{\ast}(3)$:

\begin{center}
\begin{picture}(100,70)(20,0)
\put(-15,50){$BC_{5}(3)$}
%\put(-95,45){\vector(0,-1){40}}
\put(0,45){\vector(0,-1){40}}
\put(+5,25){$\psi_{5}(3)$}
\put(-28,-5){$B_{3}(F)$}
\put(35,54){\vector(1,0){35}}
\put(70,50){$BC_{4}(3)$}
\put(85,45){\vector(0,-1){40}}
\put(90,25){$\psi_{4}(3)$}
%\put(-75,0){\vector(1,0){40}}
\put(30,0){\vector(1,0){30}}
\put(60,-5){$ B_{2}(F)\otimes F^{\ast}$}
\put(130,0){\vector(1,0){20}}
\put(110,55){\vector(1,0){40}}
\put(155,50){$BC_{3}(3)$}
\put(155,-5){$\Lambda^{3}F^{\ast}$}
\put(170,45){\vector(0,-1){40}}
\put(175,25){$\psi_{3}(3)$}
\put(260,25){(3.19)}
\end{picture}
\end{center}
\setcounter{equation}{19}
by setting it to be zero on the groups $C_{*}(k)$ for $k>3$ and using the
formulas above for the map on the subcomplex $C_{*}(3)$.

\begin{theorem} \label{z2o}
The map $\psi$ is  a morphism of
complexes.
\end{theorem}

{\bf Proof}. See chapter 3
in [G2].

{\bf 7. Complex of affine flags over a field ([G3], \S 3)}. A
p-flag in a vector
space $V$ is a sequence of subspaces
$$
0 = L^0 \subset L^1 \subset L^2\subset ...  \subset L^p \qquad dimL^i =i
$$
An affine $p$-flag is a $p$-flag $L^{\bullet}$ together with
a choice of vectors $l^i \in
L^i/L^{i-1}$ for all $1 \leq i \leq p$. We will denote affine $p$-flags as
$(l_1,...,l_p)$.

Several affine $p$-flags are in general position if all the corresponding
subspaces $L^i$ are in generic position.

Let $A^p(m)$  be the manifold of all affine $p$-flags in an $m$-dimensional
vector space $V^m$ over a field $F$. The group $GL(V^m)$ acts on it.

Let $X$ be a $G$-scheme. Then there is a simplicial scheme
$BX_{\bullet}$ where $BX_{(k)}:= G\backslash X^{k+1}$. Let $\tau_{\geq
n}BX_{\bullet}$ be the $n$-truncated simplicial scheme, where
$\tau_{\geq n}BX_{(k)} =0$ for $k <n$ and $BX_{(k)}$ otherwise.

Let $\hat BA^p(m)_{\bullet} \subset BA^p(m)_{\bullet}$ be the simplicial scheme
where $\hat BA^p(m)_{(k)}$ consists of {\it configurations} of $(n+1)$-tuples
of affine p-flags in
generic position in $V^m$ (i.e. $(n+1)$-tuples considered modulo the
action of $GL(V^m)$).

Further, let me recall  the definition of the bi-Grassmannian $\hat
G(n)$ ([G3]).
Let $(e_0,...,e_{k+l})$ be a basis in a vector space $V$. Denote by
$\hat G_l^{k}$ the open part of the Grassmannian consisting of
$l$-dimensional subspaces in $V$
transversal to the coordinate hyperplanes. It is canonically isomorphic
to the set of all $l$-planes in $k+l$-dimensional {\it affine} space $A^{k+l}$
transversal to a given $k+l$-simplex. Indeed, consider the affine hyperplane
in $V$ passing through the ends of basis vectors $e_0,...,e_{k+l}$.
There is canonical isomorphism
$m: \hat G_l^{k} \rightarrow$ $\{$ configurations of $k+l+1$ vectors
 in general position in a k-dimensional vector space $\}$.
  The  configuration $m(\xi)$ consists of  the images of $e_i$ in $V/\xi$.

The bi-Grassmannian $\hat G(n)$
 is the following diagram of manifolds
$$
\begin{array}{cccccccccccc}
&&&&&&&&&\\
&&&&&& \hat G_1^{n+2}&\stackrel{\rightarrow}{\stackrel{...}{\rightarrow}}&\hat
G_0^{n+2}\\
&&&&&&&&&\\
&&&&&&\downarrow...\downarrow&&\downarrow...\downarrow\\
&&&&&&&&&\\
\hat G(n):=&&&&\hat
G_2^{n+1}&\stackrel{\rightarrow}{\stackrel{...}{\rightarrow}}&\hat
G_1^{n+1}&\stackrel{\rightarrow}{\stackrel{...}{\rightarrow}}&\hat G_0^{n+1}\\
&&&&&&&&&\\
&&&&\downarrow...\downarrow&&\downarrow...\downarrow&&\downarrow...\downarrow\\
&&&&&&&&&\\
&&...&\stackrel{\rightarrow}{\stackrel{...}{\rightarrow}}&\hat
G_2^{n}&\stackrel{\rightarrow}{\stackrel{...}{\rightarrow}}&\hat
G_1^{n}&\stackrel{\rightarrow}{\stackrel{...}{\rightarrow}}&\hat G_0^{n}\\
\end{array}
$$
Here the horizontal arrows are provided by intersection with coordinate
hyperplanes and the vertical ones by factorisation along coordinate
axes. The bi-Grassmannian $\hat G(n)$ is a truncated simplicial scheme:
$\hat G(n)_{(k)}:= \coprod_{p+q = k} \hat G_p^{q}$.

{\bf Remark}. The bi-Grassmannian $\hat G(n)$ is not a bisimplicial
scheme. It is a {\it hypersimplicial} scheme. To explain what  it
means let me recall that the hypersimplex $\Delta^{k,l}$ is the convex
hull of centers of k-faces of the standard simplex $\Delta^{k+l+1}$
([GGL]). Its boundary is a union of hypersimplices of type $\Delta^{k-1,l}$ and
$\Delta^{k,l-1}$. More precisely, if $A$ and $S$ are 2 disjoint finite
sets, $\Delta^{A;S}$ is defined as convex hull of centers of all those
$k+|S|$-dimensional faces of $\Delta^{A\cup S}$ which contain all
vertices of $S$. Then
$$
\partial \Delta^{k,l}(e_0,...,e_{k+l}) =
\sum(-1)^i\Delta^{k-1,l}(e_0,...,\hat e_i,...,e_{k+l};e_i) +
$$
$$
\sum(-1)^i\Delta^{k,l-1}(e_0,...,\hat e_i,...,e_{k+l})
$$

{\bf Exercise}. Define hypersimplicial sets, schemes, ..., and check that
bi-Grassmannian is a $(0,n)$-truncated hypersimplicial  scheme.

A  correspondence between simplicial schemes
$X_{\bullet}$ and $Y_{\bullet}$ is a simplicial subscheme $Z_{\bullet}
\subset X_{\bullet} \times
Y_{\bullet}$ finite over $X_{\bullet}$.

There is the following   correspondence $T$ between the truncated simplicial
schemes $\tau_{\geq n}B\hat
A^p(m)_{\bullet}$ and $\hat G(n)_{\bullet}$.
For a point
$$
a = (v_0^1,...,v_0^{p+1};...;v_k^1,...,v_k^{p+1}) \in
\tau_{\geq n}B\hat A^{p+1}(n+p)_{(k)}
$$
set
$$
T(a):= \cup_{q=0}^{k-n}\cup_{i_0 +...+i_k = p-q}m^{-1}(L^{i_0}_0 \oplus
...\oplus L^{i_k}_k\vert v_0^{i_0+1},...,v_k^{i_k+1})
$$

Here $(L^{i_0}_0 \oplus
...\oplus L^{i_k}_k\vert v_0^{i_0+1},...,v_k^{i_k+1})$ is the
configuration of vectors in the space  $V^m/\oplus_{s=0}^k L_s^{i_s}$
obtained by  projection of the vectors $v_0^{i_0+1},...,v_k^{i_k+1}$ and
$m^{-1}(...)$ is the corresponding point of the appropriate Grassmannian.
\begin {theorem}  \label{at}
 $T$ is a correspondence between the truncated simplicial
schemes
$\tau_{\geq n}B\hat
A^p(m)_{\bullet}$ and $\hat G(n)_{\bullet}$
\end {theorem}
{\bf Proof}.  Follows essentially from the proof of the Key lemma in  s.2.1 of
[G3].

Let $X$ be a set. Denote by $\Bbb Z[X]$ the  free abelian group generated
by the points of $X$. Applying the functor $X \rightarrow \Bbb Z[X(F)]$ to
 our simplicial schemes  we get simplicial free  abelian
groups  $C_{\bullet}(A^p(m))$ and $BC_{\bullet}(n)$. After
normalisation we get the complex $C_{\ast}(A^p(m))$ of affine flags in
generic position
:
\begin{equation} \label{3.9m}
... \stackrel{d}{\longrightarrow} C_{n+1}(A^p(m))
\stackrel{d}{\longrightarrow} C_n(A^p(m))
\stackrel{d}{\longrightarrow} C_{n-1}(A^p(m))
\stackrel{d}{\longrightarrow} ...
\end{equation}
and the bi-Grassmannian complex $BC_{\ast}(n)$. Theorem (\ref{at})
transforms to
\begin {theorem}   \label{at1} There is a homomorphism of complexes
$$
T:C_{\ast}(A^{p+1}(p+n)) \longrightarrow BC_{\ast}(n)
$$
\end {theorem}

One has a canonical homomorphism
\begin{equation} \label{3.9b'}
H_{\ast}(GL_m(F),\Bbb Z) \longrightarrow H_{\ast}(C_{\ast}(A^{p}(m)))
\end{equation}
So we get for any $p \geq 0$ canonical homomorphisms
\begin{equation} \label{s2}
H_{\ast}(GL_{n+p}(F),\Bbb Z) \longrightarrow H_{\ast}(BC_{\ast}(n))
\end{equation}

It is sufficient  for our purposes to
consider homomorphism (\ref{s2}) for   sufficiently big $p$.

Now we need to make a statement comparing homology of the complex $BC_{*}(3)$
and cohomology of the complex $\Gamma(F,3)$. For this reason we will introduce
the cohomological version $BC^{*}(3)$ of the complex $BC_{*}(3)$ setting
$BC^{i}(3) := BC_{6-i}(3)$ and keeping the same differential, now considered as
a cohomological one. One can do the same trick with the complex
$C_{\ast}(A^{p}(m))$, getting its cohomological version $C^{\ast}(A^{p}(m))$.

Combining the map ({\ref{s2}) with
$$
\psi^{\ast}: H^{\ast}(BC_{\ast}(3)) \longrightarrow
H^{\ast}(  \Gamma(F,3))
$$
we get the desired homomorphisms
$$
H_{4}(GL(F),\Bbb Q) \longrightarrow H^2(   \Gamma(F,3) \otimes \Bbb Q)
$$

{\bf 8. The affine flag complexes over ${\cal O}$}. We will
construct in the  affine flag complex
$C_{\ast}A^{p+1}(3+p)$
a natural subcomplex $C_{\ast}A^{p+1}({\cal O},3+p)$
 corresponding to the
ring ${\cal O}$ such that

1) One has canonical homomorphism
\begin{equation}
H_{\ast}(GL_m({\cal O}),\Bbb Z) \longrightarrow H_{\ast}(C_{\ast}(A^{p}({\cal
O},m)))
\end{equation}
together with commutative diagram
$$
\begin{array}{ccc}
H_{\ast}(GL_m({\cal O}),\Bbb Z) &\longrightarrow &H_{\ast}(C_{\ast}(A^{p}({\cal
O},m)))\\
\downarrow && \downarrow\\
H_{\ast}(GL_m(K),\Bbb Z) &\longrightarrow& H_{\ast}(C_{\ast}(A^{p}(m)))
\end{array}
$$
(the down arrows are provided by the natural  embedding ${\cal O}
\hookrightarrow K$).

2) The restriction of the composition  $\psi \circ T$ to the
subcomplex \linebreak
$C_{\ast}A^{p+1}({\cal O},3+p)$ lands in  $\Gamma({\cal O},3)$:
$$
  \psi \circ T: C_{\ast}A^{p+1}({\cal O},3+p)
 \longrightarrow \Gamma({\cal O},3)
$$
In particulary this implies that the composition
$$
C_{\ast}A^{p+1}({\cal O},3+p) \longrightarrow
\Gamma(K,3) \longrightarrow \Gamma(k,2)[-1]
$$
is zero.

We will represent a $p+1$-flag in $p+3$-dimensional vector space by
vectors $(l_1,...,l_{p+1})$; the subspaces of the flag are
given by $<l_1,...,l_k>$.
Consider $m$ affine flags $a_1,...,a_m$. To define them
$(p+1) \cdot m$ vectors is  needed. We would like to define a  class
of {\it admissible} set of vectors among them. Namely
take {\it first}
$k_1$ vectors from the flag $a_1$, then first $k_2$ vectors from $a_2$ and
so on. The set of vectors we get this way is called an admissible set
of vectors related to the affine flags $a_1,...,a_m$.

Choose a basis
$e_1,...,e_{p+3}$ in ${\cal O}^{p+3}$. Let us say that  the affine flags
$a_1,...,a_m$ are in ${\cal O}$-generic
position if any admissible $p+3$-tuple of vectors  are in generic
position. This
just means that $\Delta(v_1,...,v_{p+3}) \in {\cal O}^*$ for every
admissible $p+3$-tuple of vectors $v_1,...,v_{p+3}$ related to the affine
flags $a_1,...,a_m$.

The affine flags in ${\cal O}$-generic
position provide a complex $C_{\ast}A^{p+1}({\cal O},3+p)$ with the all
described above properties.

The condition 2) holds for the following
reason. To compute $\psi \circ T$ we take a set of admissible vectors
$(v_1,...,v_{p},...)$, construct from
them  a
configuration $(v_1,...,v_{p}|v_{p+1},...)$ in a 3-dimensional
space and then apply one of homomorphisms  $f_*(3)$. The
homomorphisms  $f_*(3)$ were defined explicitely using only products
 and ratios of determinants $\Delta^{\omega}(x,y,z)$. To compute such a
determinant
we need to choose a volume form $\omega$ in the three dimensional vector
space, and
the result (homomorphisms  $f_*(3)$) does not depend on that choise. So for
each individual
configuration coming as described above from an admissible configuration
of vectors $(v_1,...,v_{p},...)$, one can choose a specific
volume form setting $\Delta^{\omega(v)}(x,y,z):=
\Delta(v_1,...,v_{p},x,y,z)$. Then for affine flags in ${\cal O}$-generic
position all the determinants we need will be in
${\cal O}^*$.

\section {Proof of Deninger's conjecture}

{\bf 1. A regulator from  $\Gamma(X,3)$ to $\Bbb R(3)_{{\cal
D}}$([G2-G3])}.  We have defined
complexes
$  \Gamma(X,3)$ so far only when $X = Spec (F)$ or $X$ is a curve.
In general $\Gamma(X,3)$  is the total complex associated
with the bicomplex
\begin {equation} \label {ga}
  \Gamma(F(X),3) \stackrel{\partial_x}{\longrightarrow} \coprod_{x
\in X_1}   \Gamma(F(x),2)[-1]
\stackrel{\partial_x}{\longrightarrow} \coprod_{x
\in X_2} F(x)^*[-2]
\stackrel{\partial_x}{\longrightarrow} \coprod_{x
\in X_3} \Bbb Q[-3]
\end {equation}

Let $S^{i}(X)$ be the space of smooth $i$-forms at the generic
point of
$X$. (This means that
each is defined on a Zariski open domain of $X$).

For any variety $X$ over $\Bbb C$ one has  {\it canonical} homomorphism of
complexes
$$
\begin{array}{ccccc}
B_{3}(\Bbb C(X))&\stackrel{\delta}{\rightarrow}& B_{2}(\Bbb C(X))\otimes
\Bbb C(X)^{\ast} & \stackrel{\delta}{\rightarrow} &
 \wedge^{3}\Bbb C(X)^{\ast}\\
&&&&\\
\downarrow r_3(1)& &\downarrow r_3(2)& &\downarrow r_3(3)\\
&&&&\\
S^0(X) &\stackrel{d}{\rightarrow}& S^1(X)
&\stackrel{d}{\rightarrow}& S^{2}(X)\\
\end{array}
$$
given by the  following formulas ($\alpha(f,g)$ was defined in (\ref{al11})).
\begin{eqnarray*}
& & r_{3}(1): \{ f\}_{3}\mapsto {\cal L}_{3}(f)\\
& & r_{3}(2) : \{ f\}_{2} \otimes g\longmapsto - {\cal
L}_{2} (f) d \arg g + \quad +\frac{1}{3} \log \vert g\vert \cdot
\alpha(1-f,f)\\
& & r_{3}(3) : f_{1}\wedge f_{2}\wedge f_{3}\mapsto {\rm Alt}
\left( \frac{1}{2} \cdot \log \vert f_{1} \vert d \arg
f_{2}\wedge d \arg f_{3} - \right. \\
& & \qquad\left. -\frac{1}{6}\log \vert f_{1} \vert
d\log \vert f_{2}\vert d\log \vert f_3 \vert \right)
\end{eqnarray*}
 It enjoys the  properties

a) $ d r_3(3)(f_1 \wedge ... \wedge f_3) + \pi_3 d\log f_1 \wedge
... \wedge d\log f_3 = 0$
where $\pi_3$ means real part.

b) Let $Y$ be an irreducible divisor in $X$ and $v_Y$ be the
corresponding valuation
on the field $\Bbb C(X)$. Then $r_3(\cdot)$  carries the   residue
homomorphism $\partial_{v_Y}$
to the usual residue homomorphism on the
DeRham complex $S^{\ast}(X) \longrightarrow S^{\ast-1}(Y)[-1]$.

A similar homomorphisms exists for the complexes $\Gamma(X,2)$ and
$\Gamma(X,1)$ (see [G2-G3] or do it as an easy exercise)

{\it This just means that these formulas provide a homomorphism from the
complex $\Gamma(X,3)$ to the weight 3 Deligne complex $\Bbb R(3)_{{\cal
D}}$ on $X$. }

{\bf 2. Relation with Beilinson's regulator}. Recall that we have
constructed in s.2-4 canonical homomorphisms
$$
c_{i,3}: K_{6-i}(\Bbb C(X))_{\Bbb Q} \longrightarrow H^i(\Gamma(\Bbb C(X)),
3)_{\Bbb Q}
$$
and in this section
$$
r_{{\cal D}}: H^i(\Gamma(\Bbb C(X)), 3)_{\Bbb Q} \longrightarrow
H^i_{{\cal D}}( Spec \Bbb C(X)), \Bbb R(3))
$$
\begin {theorem}
\label {breg}
The composition
$$
r_{{\cal D}} \circ c_{i,3}: K_{6-i}(\Bbb C(X))_{\Bbb Q} \longrightarrow
H^i_{{\cal D}}( Spec \Bbb C(X)), \Bbb R(3))
$$
coincides with Beilinson's regulator.
\end {theorem}

To prove this theorem we will remind an explicit construction of
 the universal  Chern
class $c^{\cal D}_3 \in H_{{\cal D}}^6(BGL_{\bullet}, \Bbb R(3))$ given in
[G3].
We will first construct the corresponding ``motivic'' class
$c_3 \in H_{{\cal M}}^6(BGL_{\bullet}, \Gamma(3))$ and then
apply canonical homomorphism from $\Gamma(3)$  to Deligne
cohomology.

{\bf 2. Explicit construction of the  class $c_3
\in H_{{\cal M}}^6({BGL_{3}}_{\bullet},   \Gamma(3))$ }.
  Recall that
\begin{center}
\begin{picture}(400,40)(-30,0)
\put(0,20){$BG_{\bullet}:= pt$}
\put(65,30){$\scriptstyle s_{_{0}}$}
\put(80,25){\vector(-1,0){25}}
\put(80,20){\vector(-1,0){25}}
\put(65,15){$\scriptstyle s_{_{1}}$}
\put(85,20){$G$}
\put(105,33){$\scriptstyle s_{_{0}}$}
\put(120,28){\vector(-1,0){23}}
\put(120,23){\vector(-1,0){23}}
\put(120,18){\vector(-1,0){23}}
\put(105,13){$\scriptstyle s_{_{2}}$}
\put(125,20){$G^{2}$}
\put(160,35){$\scriptstyle s_{_{0}}$}
\put(175,30){\vector(-1,0){25}}
\put(175,25){\vector(-1,0){25}}
\put(175,20){\vector(-1,0){25}}

\put(175,15){\vector(-1,0){25}}
\put(160,10){$\scriptstyle s_{_{3}}$}
\put(180,20){$G^{3}$}
\put(215,35){$\scriptstyle s_{_{0}}$}
\put(230,30){\vector(-1,0){25}}
\put(209,22){\large$\ldots$}
\put(230,15){\vector(-1,0){25}}
\put(215,10){$\scriptstyle s_{_{3}}$}
\end{picture}
\end{center}
\vskip 3mm \noindent

Choose an affine flag $a \in A^{p+1}(n+p)$. Consider simplicial
subscheme $B\hat
GL(n+p)_{\bullet} \subset BGL(n+p)_{\bullet}$ consisting of simplices
$(g_0,...,g_k)$ such that $(g_0a,...,g_ka)$ is in generic position. So
there is a morphism of
simplicial schemes
$$
A: BGL(n+p)_{\bullet} \longrightarrow B\hat A^{p+1}(n+p)_{\bullet}
$$
 defined by formula
$(g_0,...,g_k) \longrightarrow (g_0a,...,g_ka)$. Further, in s.2.? we
have constructed a morphism of truncated simplicial schemes

$$
\tau_{\geq n}\hat BA^{p+1}(n+p)_{\bullet} \longrightarrow \hat G(n)_{\bullet}
$$
So we get a morphism of truncated simplicial schemes
\begin {equation} \label{pull}
A: \tau_{\geq n}BGL(n+p)_{\bullet} \longrightarrow \hat G(n)_{\bullet}
\end {equation}
Our formulas for the homomorphism  of complexes (see s.2-4)
$$
BC^{\ast}(3)(F) \stackrel{\psi}{\longrightarrow}   \Gamma(F,3)
$$
 give us a cocycle representing a cohomology class
in $H_{{\cal M}}^6(\hat G(n)_{\bullet},   \Gamma(3))$. So pulling
 it back by (\ref{pull}) we get a cocycle $\hat c_3$ representing
$H_{{\cal M}}^6(\hat BGL(n+p)_{\bullet},   \Gamma(3))$. It is not
 a cocycle on the whole $BGL(n+p)_{\bullet}$
because it has nontrivial residues on some divisors in the complement of
 $ \hat BGL(n+p)_{\bullet}$ in $BGL(n+p)_{\bullet}$. Fortunately
it is easy to check that all residues of the components of
$\hat c_3$
in $  \Gamma(G^i,3)$ are zero for $i>3$. For $i=3$ there are
nontrivial residues, but the
corresponding problem was already solved in s.4.2 of [G3].  (Recall that by
construction components of $\hat c_3$
on $G^i$ for $i < 3$ are zero. )
Namely, in s.4 of [G3] there was  constructed a cocycle $c^M_n$
representing the
Chern class
in $H^{2n}(BGL(n+p)_{\bullet}, {\cal K}^M_n[-n])$. Here ${\cal K}^M_n$
is the sheaf
of Milnor's K-groups. In our case $n=3$ and
the component of $\hat c_3$ on $G^3$ coincides with the one of $c^M_3$.
Therefore we can simply add all components of $c^M_3$ on $G^i$
for $i <3$ and the new cochain we get will be a cocycle. Moreover,
the canonical morphism
$$
H_{{\cal M}}^6(BGL(3+p)_{\bullet},   \Gamma(3)) \longrightarrow
H^{6}(BGL(3+p)_{\bullet}, {\cal K}^M_3[-3])
$$
provided by the obvious morphism $  \Gamma(F,3) \longrightarrow
K^M_3(F)[-3]$ carries $ c_3$ to $c^M_3$ just by the construction.
In particular  the cohomology class of $ c_3$ is nonzero.
Now we apply the constructed  map to Deligne cohomology $  \Gamma(X,3)
\longrightarrow \Bbb R(3)_{{\cal D}}$ and
get a cocycle $c^{{\cal D}}_3$ representing a class in \linebreak
$H_{{\cal D}}^6(BGL(3+p)_{\bullet}, \Bbb R(3))$. It was proved in [G3]
(see s.5.7 in  [G3]) that the image of class $[c^M_3]$ in
$H^3(BGL(3+p)_{\bullet}, \Omega^3_{cl})$ coincides with the
 Chern class of universal bundle over $BG_{\bullet}$. So the
 commutative diagram
($\Omega^3_{cl} \hookrightarrow \Omega^{\geq 3}$)
$$
\begin{array}{ccc}
H_{{\cal M}}^6(BGL(3+p)_{\bullet},   \Gamma(3)) & \longrightarrow
&
H^{3}(BGL(3+p)_{\bullet}, {\cal K}^M_3)\\
&&\\
\downarrow r_{{\cal D}}&&\downarrow d\log^{\wedge 3}\\
&&\\
H_{{\cal D}}^6(BGL(3+p)_{\bullet}, \Bbb Q(3))& \longrightarrow &
H^3(BGL(3+p)_{\bullet}, \Omega^{\geq 3})
\end{array}
$$
implies
 \begin {theorem}
The cohomology class $[c^{{\cal D}}_3] \in H_{\cal
D}^6(BGL(3+p)_{\bullet}, \Bbb R(3))$ coincides with the
third Chern class of
the universal bundle.
\end {theorem}

{\bf 4. Proof of Theorem (\ref{breg})}.
 Let me first recall the definition of Beilinson's regulator
for affine  schemes.
Let $X$ be an affine  scheme over $k$ and $BGL_{\bullet}$ be the
simplicial scheme representing the classifying  space for the group $GL$.
Then $Hom_{Sch}(X,BGL_{\bullet}) = BGL(X)_{\bullet}$ is a simplicial
set. We will treat it as a 0-dimensional simplicial scheme. So
 one

has canonical morphism of simplicial schemes:
$$
X \times BGL(X)_{\bullet} \longrightarrow BGL_{\bullet}
$$
In particular we have canonical morphism
$$
e: X(\Bbb C) \times BGL(X)_{\bullet} \longrightarrow BGL_{\bullet}(\Bbb C)
$$
If $c_{n} \in H^{2n}_{{\cal D}}(BGL_{\bullet}(\Bbb C),\Bbb R(n))$ is the
universal Chern class in Deligne cohomology, then
$$
e^{\ast}c_{n}\in H^{2n}_{{\cal D}}(X(\Bbb C) \times
BGL_{\bullet}(X),\Bbb R(n))
$$
Therefore we get a homomorphism
$$
<e^{\ast}c_{n},\cdot>:H_{i}(GL(X),\Bbb Z) \longrightarrow
H^{2n-i}_{{\cal D}}(X(\Bbb C),\Bbb R(n))
$$
In particular composed with the Hurevitc map $K_i(X) \longrightarrow
H_{i}(GL(X),\Bbb Z)$ it leads to the Beilinson' regulator
$$
r_{Be}: K_i(X) \longrightarrow H^{2n-i}_{{\cal D}}(X,\Bbb R(n))
$$

Now suppose we have an $i$-cycle $\gamma$ in the complex obtained by
normalisation of the
simplicial set $BGL(\Bbb C [X])_{\bullet}$. Then to compute
$<e^{\ast}c_{n},[\gamma]> \in H^{2n-i}_{{\cal D}}(X(\Bbb C),\Bbb R(n))$
one can proceed as follows. Let $\gamma = \sum
n_j(g_0^{(j)},...,g_i^{(j)})$. Each $(g_0^{(j)},...,g_i^{(j)})$ defines a
map $\gamma_j: X(\Bbb C) \longrightarrow G^{i}$. Let $c_{{\cal D}}^i$
be the component in $\Bbb R\Gamma(G^i,\Bbb R(n)_{{\cal D}})$ of the
cocycle representing the Chern class in the
bicomplex  $\Bbb R\Gamma(BG_{\bullet},\Bbb R(n)_{{\cal D}})$.
Then $\sum_j \gamma_j^{\ast} c_{{\cal D}}^i \in \Bbb R\Gamma(X(\Bbb
C),\Bbb R(n)_{{\cal D}})$ is a cocycle representing the  class
$<e^{\ast}c_{n},[\gamma]>$.

The last problem is that the cocycle $c^{{\cal D}}_n$ is represented by
currents
on $BG_{\bullet}$, so there might be a trouble with pulling it back
by  $\gamma_j$. However on a certain generic part of
$U \subset BG_{\bullet}$ the cocycle
$c^{{\cal D}}_n$ is represented by smooth forms. We will  show that

{\it one can
always find
such a representative  $ \tilde  \gamma_j$ for the homology class class
$[\gamma_j]$ that
$\tilde \gamma_j (X(\Bbb C)) \subset U \subset G^i$}.

Currents can always be restricted to an open part of a manifold
thanks to  the map $C^{\infty}_0(U) \rightarrow C^{\infty}_0(X)$.
So presenting
$\tilde \gamma_j$ as a composition
$X(\Bbb C) \hookrightarrow U \hookrightarrow G^i
$and using  pull back of currents for open embeddings we see that
 $\gamma_j^{\ast}[c^{{\cal D}}_n]$ is represented by
$\tilde \gamma_j^{\ast}c^{{\cal D}}_n\vert_{U}$.

Now let us prove the formulated above statement.
Let $a\in V^n$, $G: = GL(V^n)$. Say that $(m+1)$-tuple of elements
$(g_0,...,g_{m+1})$ of $GL_n(F)$ is $a$-generic if the $(m+1)$-tuple of vectors
$(g_0a,...,g_{m+1}a)$ in
$V^n$ is in generic position, i.e. any $k \leq n$ of these
vectors generate a $k$-dimensional subspace.

Let $G^{m+1}(a) \in G^{m+1}$ be the subset of $a$-generic $(m+1)$-tuples of
elements. Then  $\Bbb Z[G^{m+1}(a)]$ is  a simplicial abelian group and the
corresponding complex is a free resolution of the trivial $G$-module $\Bbb Z$.
(Standard proof: if $\sum n_i(g_0^{(i)},...,g_m^{(i)})$ is a cycle, choose
an $g$ such that $ga$ is in generic position with all $g_0^{(i)}a$. Then the
boundary of
$\sum n_i(g,g_0^{(i)},...,g_m^{(i)})$ is $\sum n_i(g_0^{(i)},...,g_m^{(i)})$.
 Therefore $H_{\ast}(G,\Bbb Z) =
H_{\ast-1}( \Bbb Z [G^{\bullet}(a)]_G,\Bbb Z)$, i.e.
{\it all homology classes of $G$ can be represented by $a$-generic cycles}.
(In fact the above argument shows that this statement is true for any
reasonable notion of generic cycles. ) Theorem (\ref{breg}) is proved.

{\bf Remark}.
Similary one can define a version of continuos cohomolgy of the Lie
group $G$ as follows:
$$
H^{\ast}_{a-c}(G,\Bbb R) := H^{\ast-1}(C(G^{m+1}(a)^G)
$$
where $C(G^{m+1}(a))$ is the space of continuos functions on $G^{m+1}(a)$.
The restriction  map
$$
H_c^{\ast}(G,\Bbb R) \longrightarrow H^{\ast}_{a-c}(G,\Bbb R)
$$
is an isomorphism. Indeed, $H_c^{\ast}(G,\Bbb R) =
H_c^{\ast}(G,C(G^{\bullet}(a))$. The spectral sequence for computation of the
last group degenerates to the complex $C(G^{\bullet}(a)^G$ because
$H_c^{i}(G,C(G^{m}(a))=0$ for positive $i$ by Shapiro lemma.

Moreover the obvious pairing
$$
H_{\ast-1}( \Bbb Z [G^{\bullet}(a)]_G) \times  H^{\ast-1}(C(G^{m+1}(a)^G)
\longrightarrow \Bbb R
$$
coincides with the natural pairing $H_{\ast}(G,\Bbb Z) \times H_c^{\ast}(G,\Bbb
R) \longrightarrow \Bbb R$ after identification of the left sides.

{\bf 5. Computations for curves over $\Bbb C$}.
\begin {theorem}  \label{uhi}
  Let $X$ be a compact curve over $\Bbb C$ and
$\omega$ is a
holomorphic 1-form on $X$. and  $f_i,g_i \in
\Bbb C(X)$ are rational functions. Then
\begin{equation} \label{5.1x}
\int_{X} r_{3}(2)(\sum_i \{ f_i \}_{2} \otimes g_i) \wedge \bar \omega =
c_3 \cdot \int_{X} \log \vert g_i \vert \alpha(1-f_i,f_i) \wedge \bar \omega
\end{equation}
where $c_3 \in \Bbb Q^{\ast}$ is a constant
and  $x_i,y_i,z_i$ are divisors of functions $g_i, f_i, 1-f_i$.
\end {theorem}

{\bf Proof}. One has $F(x) d \log   g
\wedge   \omega = 0$ for a function $F(x)$ on $X$, and so
\begin{equation} \label{3.11}
 F(x) d \arg g \wedge   \omega = i F(x) \cdot  d
\log g \wedge   \omega
\end{equation}
Therefore
$$
\int_{X(\Bbb C)} {\cal L}_{2} (f )  d \arg g \wedge   \omega =i \cdot \int_{X}
{\cal L}_{2} (f ) d \log \vert g \vert \wedge   \omega
= -i \cdot \int_{X(\Bbb C)}
d {\cal L}_{2} (f ) \log \vert g \vert \wedge   \omega
$$
Here we can integrate by parts because ${\cal L}_{2} (f )$ has only
integrable singularities. Applying  the formula
\begin{equation} \label{d2}
d {\cal L}_{2} (f ) = -\log \vert 1- f\vert  d \arg \vert f\vert  +
\log \vert  f\vert d \arg \vert 1-f\vert
\end{equation}
and (\ref{3.11}) we get the proof for $n = 3$.

We did not use the crucial condition
$\sum_i(1-f_i)\wedge f_i \wedge g_i = 0 \quad \mbox{in} \quad
\Lambda^3 \Bbb Q(E)^{\ast}$ in this
computation.

{\bf 6. The case of  elliptic curves over $\Bbb C$}.
\begin {theorem}   \label{uhhh}
  Let $E$ be an elliptic curve over  $\Bbb C$ and
$\omega \in \Omega^1(\bar E)$ is normalized by $\int_{E(\Bbb C)}\omega \wedge
\bar \omega = 1$. Suppose    $f_i,g_i \in
\Bbb C(E)^{\ast}$  satisfy the  condition
\begin{equation}   \label{c3q}
\sum_i(1-f_i)\wedge f_i \wedge g_i = 0 \quad \mbox{in} \quad
\Lambda^3 \Bbb Q(E)^{\ast}
\end{equation}
Then
\begin{equation} \label{5.1q}
  \int_{E(\Bbb C)} \log \vert g_i \vert \alpha(1-f_i,f_i)  \wedge \omega =
\sum'_{\gamma_1+...+\gamma_3=0}
\frac{(x_i,\gamma_1)(y_i,\gamma_2 )(z_i,\gamma_3)(\bar\gamma_3 -
\bar\gamma_{2})}{\vert\gamma_1\vert^2 \vert\gamma_2\vert^2
\vert\gamma_3\vert^2}
\end{equation}
where   $x_i,y_i,z_i$ are divisors of functions $g_i, f_i, 1-f_i$.
\end {theorem}

{\bf Proof}. For a rational function $f(z) \in \Bbb C(E)^{\ast}$ with
$div f(z) = \sum \alpha_ix_i$ one has the Fourier expansion
\begin{equation}
\log \vert f(z)\vert = \sum_i \sum_{\gamma \in \Gamma}
\frac{\alpha_i(x_i,\gamma)}{{\vert \gamma \vert}^2} + C_f, \qquad C_f
\in \Bbb R
\end{equation}
in the sence of distributions. Indeed, $\partial \bar \partial \log\vert
f(z)\vert = \sum \alpha_i \delta_{x_i} = \sum_{\gamma \in \Gamma}
\alpha_i(x_i,\gamma)$, so $\partial \bar \partial C_f =0$, and hence $C_f$ is a
constant.

The Fourier transform carries product to the convolution and $\int_{E(\Bbb C)}$
to the functional ``value at zero''. So if we suppose all the constants
$C_f$ are zero, then we immediately get formula (\ref{5.1q}) from these
properties of the Fourier transform. In general $C_f \not = 0$. However
it turns
out  condition (\ref{c3q}),  guarantee that
  (\ref{5.1q}) is independent of  ${C_f}_i, {C_g}_i$ and
${C_{1-f}}_i$. More precisely, $f \longmapsto C_f$ is a homomorphism $\Bbb
C(X)^* \longrightarrow \Bbb R$,. We will show that (\ref{5.1q})
will not change if we replace this homomorphism by a different one.
Let us prove this statement.
In fact we will prove that
(\ref{5.1q}) written for any complex curve $X$ depends only on divisors of
$f_i, g_i, 1-f_i$.

It was shown above that
$$
\sum_i \int_{E(\Bbb C)}\log |g_i|\alpha(1 - f_i,f_i)\wedge \omega =
i \sum_i \int_{E(\Bbb C)} \log |g_i| d{\cal L}_2(f_i) \wedge \omega =
$$
$$
-i\sum_i \int_{E(\Bbb C)} {\cal L}_{2} (f )  d \log |g| \wedge \omega
$$
So the left hand side of (\ref{5.1q}) does not depend on ${C_g}_i$.

Choose a basis in $V_E := \Bbb C(E)^{\ast}\otimes \Bbb Q$. Decompose
$\sum_i(1-f_i)\wedge f_i \otimes g_i$ in this basis and collect all
the terms where a given basis element $h$ appears. We get
$\sum_i (a_i \wedge h) \otimes b_i + \sum_j (c_j \wedge d_j)\otimes h$. Let
us show that   (\ref{5.1q})
is  independent of constant $C_h$.

Indeed, $\sum_j c_j \wedge d_j = \sum_k (1-s_k)\wedge s_k$
where $\sum_k \{s_k\}_2 $ is a factor with which $h$
appears in $\sum_i \{f_i\}_2 \otimes g_i$.
Therefore
$$
\sum_i\int_{E(\Bbb C)} \alpha(c_i,d_i)\wedge \omega =
\sum_i\int_{E(\Bbb C)} \alpha(1 - s_i,s_i)\wedge \omega =
\sum_i\int_{E(\Bbb C)} d{\cal L}_{2}
(f_i )   \wedge  \omega  =0
$$
Further, thanks to condition (\ref{5.1z}) one has $(\sum_i a_i  \wedge b_i +
\sum_j c_j \wedge d_j) \wedge h = 0
$ in $\Lambda^3V_E$, so $\sum_i a_i \wedge b_i = -
\sum_j c_j \wedge d_j$ and hence
$$
\sum_i \int_{E(\Bbb C)}(\log|a_i| d\log|b_i| - \log|b_i|
d\log|a_i|)\wedge \omega = \sum_i \int_{E(\Bbb C)} \alpha(c_i,d_i)\wedge
\omega = 0
$$
On the other hand
$$
\sum_i \int_{E(\Bbb C)}(\log|a_i| d\log|b_i| + \log|b_i|
d\log|a_i|)\wedge \omega
= \sum_i \int_{E(\Bbb C)} d(\log|a_i| \log|b_i|) \wedge \omega =
0
$$
So
$$
\sum_i \int_{E(\Bbb C)}\log|a_i| d\log|b_i| \wedge \omega = \sum_i \int_{E(\Bbb
C)}\log|b_i| d\log|a_i| \wedge \omega = 0
$$
Therefore the contribution of $C_h$ is $C_h \cdot \int_{E(\Bbb
C)}\log|b_i| d\log|a_i| \wedge \omega $ and so is zero.
Theorem (\ref{uhhh}) is proved.

{\bf 7. Proof of the theorem (\ref{mth})}. Let $F$ be a number field. One has
the following
commutative diagram
$$
\begin{array}{ccccc}
K_{4}(F(X))&\stackrel{c_{2,3}}{\longrightarrow}&H^2\Gamma(Spec(F(X));3)_{\Bbb
Q} &\stackrel{r_{3}(\cdot)}{\longrightarrow}& H^2_{{\cal
D}}(Spec(F(X)\otimes \Bbb R);3)\\
&&&&\\
\downarrow \tilde \delta&&\downarrow \delta && \downarrow  res\\
&&&&\\
K_{3}(k_x)&\stackrel{ c_{1,2}}{\longrightarrow}&
H^{1}\Gamma(Spec(k_x);2)_{\Bbb Q}&\stackrel{r_{2}(\cdot)}{\longrightarrow}&
H_{{\cal
D}}^{1}\Gamma(Spec(k_x \otimes \Bbb R);2)
\end{array}
$$
where $k_x$ is also a number field (a finite extension of $F$).
Here $r_{3}(\cdot)$ and $r_{2}(\cdot)$ are  the regulator constructed
explicitely by means of the polylogarithms
We proved that  $r_{3}(\cdot) \circ c_{2,3} = r_{Be}$ where $r_{{\cal
D}}$ is the Beilinson regulator to the Deligne cohomology. Further, it is
known that $r_{2}(\cdot) \circ c_{1,2} $ coincides with the Borel
regulator $r_{Bo}$ ([G2]). So we come to the commutative diagram
$$
\begin{array}{ccc}
K_{4}(F(X))&\stackrel{r_{Be}}{\longrightarrow}& H^2_{{\cal
D}}(Spec(F(X)\otimes \Bbb R);3)\\
&&\\
\downarrow \tilde \delta && \downarrow  res\\
&&\\
K_{3}(k_x)&\stackrel{ r_{Bo}}{\longrightarrow}&
H_{{\cal
D}}^{1}\Gamma(Spec(k_x \otimes \Bbb R);2)
\end{array}
$$
The map $r_{2}(\cdot)$ is injectiv.
This follows from the injectivity of the Borel regulator and the fact that
$c_{1,2}$ is an isomorphism. So the theorem
(\ref{mth}) is proved.

 {\bf 8. The end of the proof of the main theorem}.
It follows from Beilinson's theorem on regulators of modular curves
that for a modular elliptic curve $E$ over there exists an
element $\gamma_4 \in K_4(E)$ whose regulator  (up to standard
factors) is $L(E,3)$. More precisely, there exists
a covering $X \longrightarrow E$ of $E$ by a
cetain modular curve $X$ and an element $\gamma_4' \in K_4(X)$
such that acting on $\gamma_4'$ by the transfer map $K_4(X)
\longrightarrow K_4(E)$ we get an element $\gamma_4$ with the
desired property.

The definition of the element $\gamma_4'$ and
moreover the transfer map are very implicit. So we do not get
any particular information about the element $\gamma_4$. However,
applying our results stated in theorems (\ref{uhhh}),
 (\ref{uhi}),  (\ref{breg}) to this element we get theorem (\ref{0.1}).

\section {Generalizations}

{\bf 1. The groups ${\cal R}_{n}(F)$ (see s.1.4 in [G2])}.  Let us define by
induction subgroups ${\cal R}_{n}(F)\subset \Bbb Z
[P^{1}_{F}]$, $n\geq 1$.  Set
$$
{\cal B}_{n}(F):=\Bbb Z[P^{1}_{F}]/{\cal R}_{n}(F)
$$
Put ${\cal R}_{1}(F):=(\{x\} + \{y\} - \{xy\},(x,y\in
F^{\ast})$.
Then ${\cal B}_{1}(F) = F^{\ast}$.
Consider homomorphisms
\begin{eqnarray}
& & \mbox{$\Bbb Z$}[F^*]
\stackrel{\delta_{n}}{\longrightarrow}  \left\{
\begin{array}{lll}
{\cal B}_{n-1}(F)\otimes F^{\ast} &:& n\geq 3 \\
\wedge^{2}F^{\ast} &:& n=2\end{array}\right.
\nonumber \\
& & \delta_{n}:\{x\}\mapsto \left\{ \begin{array}{lll}
\{x\}_{n-1}\otimes x &:& n\geq 3 \\
(1-x)\wedge x &:& n=2\end{array}\right.\nonumber \\
& & \delta_{n}:\{1\}\mapsto 0
\end{eqnarray}

Here $\{x\}_{n}$ is the projection of $\{x\}$ in ${\cal
B}_{n}(F)$.  Set
$
{\cal A}_{n}(F):={\rm Ker}\ \delta_{n}\; .
$
 Any element $\alpha(t) = \Sigma n_{i}\{f_{i}(t)\} \in \Bbb
Z[P^{1}_{F(t)}]$ has a specialization $\alpha(t_{0}):=\Sigma
n_{i}\{f_{i}(t_{0})\}\in \Bbb Z[F^*]$, $t_{0}\in
F^*$.  If $t_{0}$ is a zero or
pole of $f_{i}(t)$, then we put $\{f_{i}(t_{0})\}:= 0$.

\begin {definition} \label {2.10}
${\cal R}_{n}(F)$ {\it is generated by
elements} {\it   $\alpha(0)-\alpha(1)$  where $\alpha(t)$ runs
through all elements of ${\cal A}_{n}(F(t))$}.
\end {definition}

\begin {lemma}
\label {1.11}
$\delta_{n}({\cal R}_{n}(F))=0$.
\end {lemma}
See proof of lemma 1.16 in [G2]. \hfil $\Box$

Therefore we get the homomorphisms
$$
\delta: {\cal B}_{n}(F) \to \left\{ \begin{array}{lll}
{\cal B}_{n-1}(F)\otimes F^{\ast} &:& n\geq 3\\
\wedge^{2}F^{\ast} &:& n=2\end{array}\right.
$$
and finally the following complex $\Gamma (F,n)$:
$$
{\cal B}_{n}\stackrel{\delta}{\rightarrow} {\cal B}_{n-
1}\otimes F^{\ast} \stackrel{\delta}{\rightarrow} {\cal B}_{n-
2}\otimes \wedge^{2}F^{\ast} \ldots \stackrel{\delta}{\rightarrow}
{\cal B}_{2}\otimes \wedge^{n-
2}F^{\ast}\stackrel{\delta}{\rightarrow} \wedge^{n}F^{\ast}
$$
where ${\cal B}_{n}\equiv {\cal B}_{n}(F)$ placed in degree
$1$ and
$
\delta:\{x\}_{p}\otimes \wedge^{n-p}_{i=1} y_{i}\longmapsto
\delta(\{x\}_{p})\wedge \wedge^{n-p}_{i=1}y_{i}
$
has degree $+1$.

{\bf 2. The regulator to Deligne cohomology}.
Let $K$ be a  field with
a discrete valuation $v$ and the residue class $\bar
k_v$. Recall that there is
 the residue homomorphism (see ... or )
\begin{equation} \label {res}
\partial_{v} :  \Gamma (K,n)\longrightarrow   \Gamma( k_v,
n-1)[-1]
\end{equation}

Recall that $S^{i}(X)$ be the space of smooth $i$-forms at the generic point of
$X$. Set
$$
\widehat{{\cal L}_{n}} (z) = \left\{ \begin{array}{ll}
{\cal L}_{n}(z) & n :\ {\rm odd} \\
i {\cal L}_{n}(z) & n: \ {\rm even} \end{array} \right.
$$

One can show that
for $n\geq 3$
\begin{equation} \label{ss}
d \widehat{{\cal L}_{n}} (z) = \widehat{{\cal L}_{n-1}}(z)
d( i \arg z)
\end{equation}
$$
- \sum^{n-2}_{k=2} \beta_{k}
\log^{k-1}\vert z \vert \cdot \widehat{{\cal L}_{n-k}}(z)
\cdot d\log \vert z\vert\nonumber \\
+ \beta_{n-1} \log^{n-2}\vert z \vert
\alpha(1-z,z)\; .
$$
In this formula the same coefficients
appear as in the definition of the function ${\cal L}_{n}$.

\begin {theorem}
\label {1.11a}
There exist {\bf canonical} homomorphism of complexes
$$
\begin{array}{ccccccc}
{\cal B}_{n}(\Bbb C(X))&\stackrel{\delta}{\rightarrow}& {\cal B}_{n-
1}(\Bbb C(X))\otimes \Bbb C(X)^{\ast} & \stackrel{\delta}{\rightarrow} & \ldots
&\stackrel{\delta}{\rightarrow} &
 \wedge^{n}\Bbb C(X)^{\ast}\\
&&&&&&\\
\downarrow r_n(1)& &\downarrow r_n(2)& & & &\downarrow r_n(n)\\
&&&&&&\\
S^0(X) &\stackrel{d}{\rightarrow}& S^1(X) &\stackrel{d}{\rightarrow}& ...
&\stackrel{d}{\rightarrow}& S^{n-1}(X)\\
\end{array}
$$
with the following properties:

a) $ d r_n(n)(f_1 \wedge ... \wedge f_n) + \pi_n d\log f_1 \wedge
... \wedge d\log f_n
= 0$
where $\pi_n  $ means real part for $n$ odd and imaginary  for $n$ even.

b) $r_n(1)\{f(x)\}_n = {\cal L}_n(f(x))$ and
\begin{equation} \label{ss1}
   r_{n+1}(n)(\{ f\}_{n-1}\otimes g) :=
\widehat{{\cal L}_{n-1}}
(f) di \arg g -
\end{equation}
$$
 -\sum^{n-2}_{k=2} \beta_k
\log^{k-2}\vert f\vert \log \vert g\vert\cdot \widehat{\cal L}_{n-k} (f)
d\log \vert f\vert
 + \beta_{n-1} \log \vert g\vert \cdot
\log^{n-3} \vert f\vert \cdot \alpha(1-f,f)\;
$$

d) Let $Y$ be an irreducible divisor in $X$ and $v_Y$ be the
corresponding valuation
on the field $\Bbb C(X)$. Then $r_n(\cdot)$  carries the   residue
homomorphism $\partial_{v_Y}$
(see (\ref{res}))  to
the residue homomorphism on the
DeRham complex $S^{\ast}(X) \longrightarrow S^{\ast-1}(Y)[-1]$.
\end {theorem}

An explicit construction of this homomorphism
will be given
elsewhere.

{\bf Remark}. It was conjectured in [G2] that the complex ${\cal
  B}_{*}(\Bbb C(X))$ computes the weight $n$ pieces of the
  $K$-theory of the field $\Bbb C(X)$. The homomorphism $r_n(\cdot)$
  should provide the regulator map to Deligne cohomology.

{\bf 3. Computations for curves over $\Bbb C$}.
\begin {theorem}
 Let $X$ be a compact curve over $\Bbb C$ and
$\omega$ is a
holomorphic 1-form on $X$. Suppose that $n > 3$ and  $f_i,g_i \in
\Bbb C(X)$ are rational functions satisfying the following condition:
\begin{equation} \label{5.1zz}
 \sum_i\{f_i\}_{n-2}\otimes f_i \wedge g_i =0 \quad \mbox{in}\quad
{\cal B}_{n-2}(\Bbb C(X))\otimes  \Lambda^2 \Bbb C(X)^{\ast}
\end{equation}
Then
\begin{equation} \label{5.1x}
\int_{X} r_{n}(2)(\sum_i \{ f_i \}_{n-1} \otimes g_i) \wedge \bar \omega =
c_n \cdot \int_{X} \log \vert g_i \vert \log
\vert f_i \vert^{n-3}  \alpha(1-f_i,f_i) \wedge \bar \omega
\end{equation}
where $c_n \in \Bbb Q^{\ast}$ is a constant
and  $x_i,y_i,z_i$ are divisors of functions $g_i, f_i, 1-f_i$.
\end {theorem}
\vskip 3mm \noindent
{\bf Proof}. Let us first give detailed proof in  the cases $n=4$
emphasizing certain differences between this case and $n=3$ case considered
in theorem (\ref{uhi}).

i) $n=4$.
$$
\int_{X(\Bbb C)} {\cal L}_{3} (f )  d \arg g \wedge   \omega = i \cdot
\int_{X(\Bbb C)}
{\cal L}_{3} (f ) d \log \vert g \vert \wedge   \omega
= -i \cdot \int_{X(\Bbb C)}
d {\cal L}_{3} (f ) \log \vert g \vert \wedge   \omega
$$
Applying  the formula for $d {\cal L}_{3} (f )$ we get
$$
 -i\cdot \int_{X(\Bbb C)}   \Bigl(  -{\cal
L}_{2} (f ) d \arg f    +\frac{1}{3} \log \vert f\vert \cdot \alpha(1-f,f)
\Bigr) \cdot  \log \vert g \vert   \wedge    \omega
$$
It seems that in general it is impossible to rewrite the individual integral
$\int_{X(\Bbb C)} {\cal L}_{2} (f ) log \vert g \vert d log \vert f \vert
\wedge    \omega$
as
$$
c \cdot\sum_i \int_{X(\Bbb C)} \log \vert t_i \vert \log
\vert s_i \vert  \alpha(1-s_i,s_i) \wedge   \omega
$$
for some rational functions $s_i$ and $t_i$.
However  assuming  condition (\ref{5.1zz})   one can do this for
\begin{equation} \label{3.a}
\sum_i\int_{X(\Bbb C)} {\cal L}_{2} (f_i ) \log \vert g_i \vert d \log \vert
f_i
\vert \wedge   \omega
\end{equation}
Indeed,  (\ref{5.1zz}) just means that  one has
$$
\sum_i \{f_i\}_2 \otimes f_i \otimes g_i   \in B_2(\Bbb C(X)) \otimes S^2
\Bbb C(X)^{\ast}
$$
The 2 homomorphisms  from
$B_2(\Bbb C(X)) \otimes \Bbb C(X)^{\ast} \otimes \Bbb C(X)^{\ast}$
to real $C^{\infty}$ 1-forms on open part of $X$ given by the formulas
$$
\{f\}_2 \otimes f \otimes g \longmapsto  {\cal L}_{2} (f ) log \vert f \vert d
log \vert g \vert
$$
and
$$
   \{f\}_2 \otimes f \otimes g \longmapsto \frac{1}{2} {\cal
L}_{2} (f ) \Bigl( \log \vert f \vert d \log \vert g \vert \quad + \quad
\log \vert g \vert
d\ log \vert f \vert  \Bigl)
$$
coincide on the subgroup $B_2(\Bbb C(X)) \otimes S^2
\Bbb C(X)^{\ast}$. Therefore (\ref{3.a}) is equal to
$$
\frac{1}{2}\sum_i\int_{X(\Bbb C)} {\cal L}_{2} (f_i )   d \Bigl(log \vert g_i
\vert
\cdot
log \vert f_i \vert \Bigr)\wedge   \omega = - \frac{1}{2}
 \sum_i\int_{X(\Bbb C)} d {\cal L}_{2} (f_i )log
\vert g_i \vert \cdot log \vert f_i \vert \wedge   \omega
$$
It remaines to use the formula (\ref{d2})
together with
(\ref{3.11}).
The theorem for $n=4$  is proved.

ii)
The  proof of the general statement is based on the following
\vskip 3mm \noindent
\begin {lemma}  Let us suppose (\ref{5.1zz}). Then
\begin{equation} \label{5.1c}
\sum_i\{f_i\}_{n-k}\otimes \underbrace{f_i\otimes ...
\otimes f_i}_{k
\quad \mbox{times}} \otimes g_i =0 \quad \in
{\cal B}_{n-k}(\Bbb C(X))\otimes  S^{k+1} \Bbb C(X)^{\ast}
\end{equation}
\end {lemma}
{\bf Proof}. Indeed, according to (\ref{5.1zz})
$$
\sum_i\{f_i\}_{n-2}\otimes f_i \otimes g_i \in {\cal B}_{n-2} \otimes
S^{2} \Bbb C(X)^{\ast}
$$
and from the other hand
$$
\sum_i\{f_i\}_{n-k}\otimes \underbrace{f_i\otimes ...  \otimes f_i}_{k
\quad \mbox{times}} \otimes g_i \in {\cal B}_{n-k-1} \otimes
S^{k} \Bbb C(X)^{\ast} \otimes \Bbb C(X)^{\ast}
$$
It remains to use the fact that for a vector space $V$
$$
S^kV \otimes V \cap S^{k-1}V \otimes S^2V = S^{k+1}V
$$
For elements $f_i,g_i$ satisfying (\ref{5.1zz}) one has
\begin{equation}   \label{masha}
\sum_i \int_{X(\Bbb C)}{\cal
L}_{n-k}(f_i)\log^{k-1}|f_i|\log|g_i|d\log|f_i|\wedge \omega  =
\end{equation}
$$
\frac{1}{k+1} \int_{X(\Bbb C)}{\cal
L}_{n-k}(f_i)d(\log^{k-1}|f_i|\log|g_i|)\wedge \omega
$$
Integrating by parts and using the formula for $d{\cal L}_{n-k}(f)$  and
(\ref{3.11}) we get the theorem by induction.

Similar arguments prove the following
\begin {proposition} \label{ma1}
If $f_i,g_i$ satisfy (\ref{5.1zz}) then for any $n-2 \geq k \geq 1$ one
has $(q_k \in \Bbb C^{\ast}$
$$
\sum_i \int_{X(\Bbb C)}\alpha(1-f_i, f_i)\log^{n-3}|f_i|\log|g_i|\wedge
\omega =
$$
$$q_k \cdot\sum_i \int_{X(\Bbb C)} d{\cal
L}_{n-k}(f_i)\log^{k-1}|f_i|\log|g_i|\wedge \omega
$$
\end{proposition}

{\bf 4. The case of  elliptic curves over $\Bbb C$}.
\begin {theorem}
  Let $E$ be an elliptic curve over  $\Bbb C$ and
$\omega \in \Omega^1(\bar E)$ is normalized by $\int_{E(\Bbb C)}\omega \wedge
\bar \omega = 1$. Suppose    $f_i,g_i \in
\Bbb C(E)^{\ast}$  satisfies the  condition
\begin{equation} \label{5.1z}
 \sum_i\{f_i\}_{n-2}\otimes f_i \wedge g_i =0 \quad \mbox{in}\quad
{\cal B}_{n-2}(\Bbb C(X))\otimes  \Lambda^2 \Bbb C(X)^{\ast}
\end{equation}
Then
\begin{equation} \label{5.1qw}
  \int_{E(\Bbb C)} \log \vert g_i \vert \log
\vert f_i \vert^{n-3} \alpha(1-f_i,f_i)  \wedge \bar \omega =
\end{equation}
$$
 \sum'_{\gamma_1+...+\gamma_n=0}
\frac{(x_i,\gamma_1)(y_i,\gamma_2 +
...+\gamma_{n-1})(z_i,\gamma_n)(\bar\gamma_n -
\bar\gamma_{n-1})}{\vert\gamma_1\vert^2 \vert\gamma_2\vert^2 ...
\vert\gamma_n\vert^2}
$$
where   $x_i,y_i,z_i$ are divisors of functions $g_i, f_i, 1-f_i$.
\end {theorem}

{\bf Proof}. It is similar to the case $n=3$ considered in theorem
(\ref{uhhh}). Recall the Fourier expansion
\begin{equation}
\log \vert f(z)\vert = \sum_i \sum_{\gamma \in \Gamma}
\frac{\alpha_i(x_i,\gamma)}{{\vert \gamma \vert}^2} + C_f, \qquad C_f
\in \Bbb R
\end{equation}

As before, assuming all the constants
$C_f$ are zero we immediately get formula (\ref{5.1qw}) from the
properties of the Fourier transform. In general $C_f \not = 0$. However
the condition (\ref{5.1z}) guarantee that
  (\ref{5.1qw}) is independent of  ${C_f}_i, {C_g}_i$ and ${C_{1-f}}_i$.
Let us prove this statement.

We will consider separately cases $n=4$ and $n > 4$
to emphasize the main points of the calculation. In fact we will prove that
(\ref{5.1qw}) written for any complex curve $X$ depends only on divisors of
$f_i, g_i, 1-f_i$.

i) $n=4$.  Consider the expression
\begin {equation} \label{hu}
\sum_i(1-f_i)\wedge f_i \otimes f_i \otimes g_i \in \otimes^4 V_E
\end {equation}
Recall that we suppose
\begin {equation} \label{huz}
\sum_i(1-f_i) \wedge f_i \otimes f_i \wedge g_i  = 0 \quad \mbox{in}
\quad \wedge^2
V_E \otimes \wedge^2 V_E
\end {equation}
Let
$s = \prod_i h_i^{<s,h_i>}$
be a notation for decomposition of a function $s$ in chosen basis $h_i$.
Then the component of (\ref{hu}) in $\otimes^3 V_E \otimes h$ is
$$
\sum_i <g_i,h> \cdot (1-f_i)\wedge f_i \otimes f_i \otimes h
$$
So the
contribution of $C_h$ is
$$
C_h \cdot \sum_i <g_i,h> \cdot \int_{E(\Bbb C)} \alpha(1-f_i, f_i) \log|f_i|
\wedge
\omega
$$
This is zero because  using proposition(\ref{ma1}) and Stokes formula
$$
\int_{E(\Bbb C)} \alpha(1-f_i, f_i) \log|f_i| \wedge
\omega = 3/2\cdot \int_{E(\Bbb C)} d{\cal L}_3(f_i) \wedge \omega=0
$$

Now consider the component of (\ref{hu}) in $\wedge^2 V_E \otimes
h\otimes V_E$. Let us write it as
$\sum_i a_i \wedge b_i \otimes h \wedge c_i$. Then condition (\ref{huz})
implies that
$$
\sum_i a_i \wedge b_i \otimes  c_i = \sum_i <g_i,h> \cdot (1-f_i)\wedge f_i
\otimes f_i
$$
Therefore
$$
\int_{E(\Bbb C)}\alpha(a_i,b_i) \log|c_i| \wedge \omega = \sum_i <g_i,h> \cdot
\int_{E(\Bbb C)} \alpha(1-f_i, f_i) \log|f_i| \wedge \omega = 0
$$

Finally,  look at the component of (\ref{hu})   in $h \wedge V_E \otimes
\otimes^2 V_E$. It actually belongs to $h \otimes V_E \otimes S^2 V_E$
because of the condition (\ref{huz}).
Let us decompose it on 2 components: the first in $h \otimes S^3 V_E$ and the
second in $h \otimes
\Bigl(\wedge^2 V_E \otimes V_E \cap V_E \otimes S^2 V_E\Bigr)$.

If we write the first component as $\sum_i h \otimes x_i \cdot y_i \cdot
z_i$, the corresponding contribution of $C_h$ will be
$$ C_h 1/3 \sum_i \int_{E(\Bbb C)} d(\log|x_i| \log|y_i| \log|z_i|) = 0
$$
It remains  the second component.The reason the
contribution of $C_h$ to be  zero in this case is the most funny.
Namely, this component  can be written as
$\sum_i h \wedge (1-s_i) \otimes s_i \otimes g_i$. So the corresponding
integral is
$$
- C_h \cdot \int_{E(\Bbb C)} d\log|1-s_i| \log|s_i| \log|g_i| \wedge
\omega =-1/2 C_h \int_{E(\Bbb C)} \alpha(1-s_i, s_i) \log|s_i| \wedge \omega
$$
(We used the fact that
  $\int_{E(\Bbb C)}d \log |a_i| \log |b_i| \log |c_i| \wedge
\omega = 1/2\int_{E(\Bbb C)} \alpha(a_i, b_i)\log |c_i|\wedge
\omega$ if $\sum_i a_i \otimes b_i \otimes c_i \Lambda^2V_E \otimes V_E$).
But this integral coincides with the one for $\sum_i(1-s_i) \wedge s_i
\otimes h
\otimes s_i$ which already was proved to be zero!

ii) $n>4$. The reasons are similar to those of the case $n=4$.
Proposition (\ref{ma1}) for $k=2$ implies the statement about $C_{g_i}$.
Consider element
\begin {equation} \label{huma}
\sum_i (1-f_i)\wedge f_i\underbrace{\otimes f_i \otimes ... \otimes
f_i}_{n-3 \quad \mbox{times}} \otimes g_i \in \Lambda^2 V_E \otimes S^{n-3}
V_E \otimes V_E
\end{equation}
The condition that its projection to $\Lambda^2 V_E \otimes \otimes^{n-4} V_E
\otimes\Lambda^2 V_E$ is zero implies that the contribution of $C_h$ related to
the term
$$
\sum_i (1-f_i)\wedge f_i\underbrace{\otimes f_i \otimes ... \otimes
f_i}_{n-4 \quad \mbox{times}} \otimes h \otimes g_i
$$
is zero (the arguments  are in complete analogy with the $n=4$ case).

Finally,   the component of (\ref{huma}) in $h \wedge V_E \otimes
\otimes^{n-2} V_E$belongs to $h \otimes V_E \otimes S^{n-2} V_E$
thanks to condition (\ref{huz}).
Decomposing  it on 2 components:  in $h \otimes S^{n-1} V_E$ and the
 in $h \otimes
\Bigl(\wedge^2 V_E \otimes^{n-3} V_E \cap V_E \otimes S^{n-2} V_E\Bigr)$.
we get the statement similarly to the case $n=4$.
Theorem is proved.

\section{Appendix}
{\bf 1. Proof of theorem \ref{z2}b)}. Let me remind the formulation of
this theorem

{\bf Theorem \ref{z2}}
a) $f_{4}(3)$  and $f_{5}(3)$ {\it do not depend on the choice
of $\omega$.

b) The homomorphisms $f_*(3)$ provide a morphism of
complexes. }

Proof. a) See the proof of  similar results in chapter 3
of [G2].

b) We have to prove that $f_4(3) \circ d = \delta \circ f_5(3)$ and
$f_5(3) \circ d = \delta \circ f_6(3)$. For the first result see chapter 3 in
[G2].

The second one is much more subtle. As pointed out H.Gangl, the geometric
proof given in  [G2] (see theorem 3.10 there) has some errors.
Namely, in lemma 3.8
$r = -r_3$ but not $r=r_3$ as clamed, and as a result the proof
of theorem 3.10 become more involved; further,  the correct
statement in theorem 3.10 is $f_5(3) \circ d =  \delta \circ 1/15 \cdot
f_6(3)$
(the coefficient $1/15$ in the definition of $f_6(3)$ was missed).

 Another proof was
given in [G1]. It was actually the first proof of the statement b).
However in this proof we
used a different definition for homomorphism $f_6(3)$
(the map $M_3$ in [G1]). Moreover the proof was rather complicated and the
relation between  the homomorphisms $f_6(3)$ and $M_3$   not easy to see.
Therefore I will present in detail a completely different proof togerther
with some corrections to chapter 3 in [G2].

Let us suppose that in a three dimensional vector space $V_3$ we choose a
volume form $\omega$. Then for any two vectors $a,b$
one can define the cross product $a \times b \in V_3^*$ as follows:
$<a \times b, c>: = \Delta(a,b,c)$. The volume form $\omega$ defines the
dual volume form in $V_3^*$, so we can define $\Delta(x,y,z)$ for
any three vectors in $V_3^*$.
\begin{lemma} \label{gz2}
For any $6$ vectors in generic position $a_1,a_2,a_3,b_1,b_2,b_3$ in $V_3$
$$
\Delta(a_1,a_2,b_1) \cdot \Delta(a_2,a_3,b_2) \cdot \Delta(a_3,a_1,b_3)  -
\Delta(a_1,a_2,b_2) \cdot \Delta(a_2,a_3,b_3) \cdot \Delta(a_3,a_1,b_1) =
$$
$$
\Delta(a_1,a_2,a_3) \cdot \Delta(a_1 \times b_1,a_2 \times b_2,a_3 \times b_3)
$$
\end{lemma}

{\bf Proof}. The left hand side is zero if the
vectors $a_1,a_2,a_3$ are linearly dependent. So $\Delta(a_1,a_2,a_3)$
divides it. Similarly the left hand side is zero if $a_i$ is
collinear to $b_i$ or $\alpha_1 a_1 + \beta_1 b_1 = \alpha_2 a_2 + \beta_2 b_2=
\alpha_3 a_3 + \beta_3 b_3$ for some numbers $\alpha_k, \beta_k$. This
implies that $\Delta(a_1 \times b_1,a_2 \times b_2,a_3 \times b_3)$ also
divides the left hand side. It is easy to deduce the formula from this.

However it perhaps easier to check the formula directly.
Consider the following special configuration of
vectors:
$$
\begin{array} {cccccc}
a_1&a_2&a_3&b_1&b_2&b_3\\
-&-&-&-&-&-\\
1&0&0&x_1&y_1&z_1\\
0&1&0&x_2&y_2&z_2\\
0&0&1&x_3&y_3&z_3
\end{array}
$$
Then the left hand side is equal to
$x_3 y_1 z_2 - y_3 z_1 x_2$, and the computation of the right hand side
gives the same result. The lemma is proved.

{\bf Remark}. Let $a_1,...,a_n,b_1,...,b_n$ be a configuration of $2n$
vectors in
an $n$-dimensional vector space $V_n$. Set $\Delta({\hat a}_n,b_1):=
\Delta(a_1,...,a_{n-1},b_1)$ and so on. Then
$$
\Delta({\hat a}_1,b_1) \cdot ... \cdot  \Delta({\hat a}_{n},b_n)  -
\Delta({\hat a}_1,b_2) \cdot ... \cdot  \Delta({\hat a}_n,b_1) =
$$
$$
\Delta(a_1,...,a_n) \cdot \Delta(a_1 \times ... \times a_{n-2} \times b_n,
... ,a_n \times ... \times a_{n-3} \times b_{n-1})
$$

Notice that $f_5(3) \circ d - \delta \circ f_6(3) \in B_2(F) \otimes F^*$.
There is a homomorphism
$$
\delta \otimes id: B_2(F) \otimes F^* \longrightarrow \wedge^2F^*\otimes F^*,
\qquad \{x\}_2
 \otimes y \longmapsto (1-x) \wedge x \otimes  y
$$
The crucial step of the proof is the following
\begin{proposition} \label{gz3}
$$
(\delta \otimes id) \circ \Bigl
(f_5(3) \circ d - \delta \circ f_6(3)\Bigr)(v_1,...,v_6) =0
\quad \mbox{in} \quad \wedge^2 F^* \otimes F^*
$$
\end{proposition}

{\bf Proof}. We will use notation
$\Delta ( i,j,k) $ for $\Delta ( v_i,
v_j, v_k)$.  According to lemma (\ref{gz2})
$$
1 - \frac{\Delta ( 1, 2, 4)\Delta (2,3, 5)
\Delta(3,1,6)}
{\Delta(1,2,5 )\Delta (2,3,6)\Delta
(3,1,4)} =
\frac{\Delta (1,2,3) \Delta (v_1 \times v_4, v_2\times v_5,
v_3\times v_6)}
{\Delta(1,2,5 )\Delta (2,3,6)\Delta
(3,1,4)}
$$
Using the cyclic permutation
$1-> 2->3 ->1, 4->5->6->4$ we see that one has to calculate  the element
$$
3 \cdot {\rm Alt}_{6}\left\{ \frac{\Delta (1,2,4)
\Delta (2,3, 5)
\Delta(3,1,6)}
{\Delta(1,2,5 )\Delta (2,3,6)\Delta
(3,1,4)} \wedge \frac{\Delta (1,2,3) \Delta (v_1 \times v_4, v_2\times v_5,
v_3\times v_6)}
{\Delta(1,2,5 )\Delta (2,3,6)\Delta
(3,1,4)} \otimes \frac{\Delta(1,2,4)}{\Delta(1,2,5)}\right\}
$$
in $\wedge^2 F^* \otimes F^*$.

Let us do this. We will compute first
the contribution of the factor $\otimes \Delta(1,2,4)$.
What we need to find is
$$
{\rm Alt}_{(1,2,4);(3,5,6)}\left\{\frac{\Delta (1,2,4)\Delta (2,3, 5)
\Delta(3,1,6)}
{\Delta(1,2, 5)\Delta (2,3,6)\Delta
(3,1,4)} \wedge \frac{\Delta (1,2,3)
\Delta (v_1 \times v_4, v_2\times v_5,
v_3\times v_6)}
{\Delta(1,2,5 )\Delta (2,3,6)\Delta
(3,1,4)}\right\}
$$
in $\wedge^2 F^*$.
Here ${\rm Alt}_{(1,2,4);(3,5,6)}$
is the skewsymmetrization with respect to the
 group $S_3 \times S_3$ which permutes the indices $(1,2,4)$ and $(3,5,6)$.

i) Consider
$$
{\rm Alt}_{(1,2,4);(3,5,6)}\left\{\frac{\Delta (1,2,4)\Delta (2,3, 5)
\Delta(3,1,6)} {\Delta(1,2,5 )\Delta (2,3,6)\Delta
(3,1,4)} \wedge \Delta (v_1 \times v_4,v_2\times v_5,
v_3\times v_6)\right\}
$$
Using the skewsymmetry with respect to the permutation  exchanging
$1$ with $3$  as well as $4$ with $ 6$  (notation: $: (13)(46)$) we see that
this expression  is zero.

ii)Look at
$$
- {\rm Alt}_{(1,2,4);(3,5,6)}\left\{\frac{\Delta (1,2,4)\Delta (2,3,5)
\Delta(3,1,6)} {\Delta(1,2,5)\Delta (2,3,6)\Delta
(3,1,4)} \wedge \Delta (2,3,6) \otimes \Delta(1,2,4)\right\}
$$
The skewsymmetry with respect to $(14)$  or with respect to
$(36)$ imply that it is also zero.

iii) Consider
$$
- {\rm Alt}_{(1,2,4);(3,5,6)}\left\{\frac{\Delta (1,2,4)\Delta (2,3, 5)
\Delta(3,1,6)} {\Delta(1,2, 5)\Delta (2,3,6)\Delta
(3,1,4)} \wedge \Delta (1,2,5) \otimes \Delta(1,2,4)\right\}
$$
The skewsymmetry with respect to the permutations $(12)$ as well as $(36)$
leads to
$$
-{\rm Alt}_{(1,2,4);(3,5,6)}\left\{\frac{\Delta(2,3,5)}{\Delta(1,3,4)}\wedge
\Delta(1,2,5)
\otimes \Delta(1 ,2,4)\right\}
$$

iv) Look at the term with $\Delta(3,1,4)$:
$$
- {\rm Alt}_{(1,2,4);(3,5,6)}\left\{\frac{\Delta (1,2,4)\Delta (2,3, 5)
\Delta(3,1,6)} {\Delta(1,2, 5)\Delta (2,3,6)\Delta
(3,1,4)} \wedge \Delta (3,1,4) \otimes \Delta(1,2,4)\right\}
$$
Using the permutation $(14)$ we get
$$
-{\rm Alt}_{(1,2,4);(3,5,6)}\left\{\frac{\Delta(3,1,6)}{\Delta(1,2,5)}\wedge
\Delta(1,3,4)
\otimes \Delta(1 ,2,4)\right\}
$$
v) Finally, using $(12)$ and  $(56)$
we see that
$$
{\rm Alt}_{(1,2,4);(3,5,6)}\left\{\frac{\Delta (1,2,4)\Delta (2,3, 5)
\Delta(3,1,6)} {\Delta(1,2, 5)\Delta (2,3,6)\Delta
(3,1,4)} \wedge \Delta (1,2,3) \otimes \Delta(1,2,4)\right\} =
$$
$$
 - 3 \cdot {\rm Alt}_{(1,2,4);(3,5,6)}
\left\{\Delta(1,3,5) \wedge \Delta(1,2,3)
\otimes \Delta(1 ,2,4)\right\}
$$

Therefore we get
$$
{\rm Alt}_{(1,2,4);(3,5,6)}\Bigl(
\Delta(1,2,5) \wedge \frac{\Delta(2,3,5)}{\Delta(1,3,4)} + \Delta(1,3,4) \wedge
\frac{\Delta(1,3,6)}{\Delta(1,2,5)} +
$$
$$
3 \cdot \Delta(1,2,3) \wedge \Delta(1,3,5)
\Bigr)\otimes \Delta(1,2,4)  =
$$
$$
{\rm Alt}_{(1,2,4);(3,5,6)} \Bigl(\Delta(1,2,5) \wedge \Delta(2,3,5) +
\Delta(1,3,4) \wedge \Delta(1,3,6) +
$$
$$
3 \cdot \Delta(1,2,3) \wedge \Delta(1,3,5)
\Bigr)\otimes \Delta(1,2,4)  =
$$
$$
5 \cdot {\rm Alt}_{(1,2,4);(3,5,6)} 1759 \Delta(1,2,3) \wedge
\Delta(1,3,5) \otimes \Delta(1,2,4)
$$
The computation of the contribution of  $\Delta(1,2,5)$ goes
similarly and gives the same answer.
{\it So the total result of our computation is}
\begin{equation} \label{su2}
-30 \cdot {\rm Alt}_{6} \left\{ \Delta(1,2,4) \wedge \Delta(1,4,5)
\otimes \Delta(1,2,3) \right\}
\end{equation}
Here we get the coefficient $-30$  taking into account the action of the cyclic
group of order $3$ generated by $1->2->3->1, 4->5->6->3$.

Now let us compute $f_5(3) \circ d(v_1,...,v_6)$. We will use the formula
\begin{equation} \label{su1}
\delta \{r(v_1,v_2,v_3,v_4)\}_2 = 1/2\cdot {\rm Alt}_{4}\left\{\Delta(v_1,v_2)
\wedge \Delta(v_1,v_3)\right\}
\end{equation}
Neglecting for a moment
the  constant $c, c'$ we get
$$
(\delta \otimes id) \Bigl(f_5(3) \circ d(v_1,...,v_6)\Bigr) =
c \cdot {\rm Alt}_{6} \{r(v_1|v_2,v_3,v_4,v_5\}_2 \otimes \Delta(1,2,3) =
$$
$$
c' \cdot {\rm Alt}_{6} \Delta(1,2,4) \wedge \Delta(1,4,5)\otimes \Delta(1,2,3)
$$
To justify this we  used here formula (\ref{su1}) and the symmetry
considerations
for transpositions $i<->j$ where $1 \leq i<j\leq 3$. More careful
consideration
shows $c' =-2$. It remains to compare it with (\ref{su2}).
That's why we need in the definition of $f_6(3)$
the coefficient $1/15$.

We have proved that
$$
(f_5(3) \circ d - \delta \circ f_6(3)\Bigr)(v_1,...,v_6) =
\sum_{1 \leq i<j<k \leq 6 } \gamma_{i,j,k} \otimes \Delta(i,j,k)
$$
where $\gamma_{i,j,k} \in B_2(F)$ and moreover $\delta(\gamma_{i,j,k}) =0$ in
$\wedge^2F^*$. According to [S2]
\begin{equation} \label{susl}
Ker\Bigl( B_2(F)
\stackrel{\delta}{\longrightarrow} \wedge^2F^*\Bigr)\otimes \Bbb Q =
K_3^{ind}(F)\otimes \Bbb Q
\end{equation}
One knows that $K_3^{ind}(F(t)) \otimes \Bbb Q = K_3^{ind}(F)\otimes \Bbb Q$.
Therefore the left hand side of (\ref{susl}) is rationaly invariant.
On the other hand one can connect  by a
rational curve the configurations
$(v_1,v_2,...,v_6)$ and $(v_2,v_1,...,v_6)$ (interchanging $v_1$ with $v_2$)
in the space of all generic  configurations.
This implies that $\gamma(1,2,3) =  \gamma(2,1,3)$ modulo torsion.
But
$\gamma(1,2,3) =  - \gamma(2,1,3)$ modulo torsion by the skewsymmetry.
So $\gamma(1,2,3) =0$ modulo torsion, and  the same conclusion is valid
for  $\gamma(i,j,k)$. With more work one can show that
$f_5(3) \circ d - \delta \circ f_6(3) = 0$ at least modulo 6-torsion,
but we do not need this. Theorem is proved.

{\bf 2. The geometrical definition of the homomorphism $f_6(3)$}  Let
$(a_{1},a_{2},a_{3},b_{1},b_{2},b_{3})$ be a configuration of
6 distinct points in $P^{2}$ as on
fig.\ 1.  Let $P^{2}=P(V_{3})$.  Choose vectors in $V_{3}$
such that they are projected to points $a_{i},b_{i}$.
We  denote them by the same letters.
Choose $f_{i} \in V_{3}^{\ast}$ such that $f_{i}(a_{i}) =
f_{i}(a_{i+1}) = 0$.  Put
\begin{equation}
r'_{3}(a_{1},a_{2},a_{3},b_{1},b_{2},b_{3}) =
\frac{f_{1}(b_{2}) \cdot f_{2}(b_{3})\cdot f_{3}(b_{1})}
{f_{1}(b_{3})\cdot f_{2}(b_{1})\cdot f_{3}(b_{2})}\; .
\end{equation}
The right-hand side of (3.10) does not depend on the choice of
vectors $f_{i},b_{j}$.
\begin{center}
\begin{picture}(100,80)
\put(47,78){$a_{2}$}
\put(47,68){$\bullet$}
\put(50,70){\vector(2,-3){40}}
\put(50,70){\vector(-2,-3){40}}
\put(23,31){$\bullet$}%
\put(12,33){$b_{1}$}
%\put(71,31){$\bullet$}
\put(91,10){\vector(-1,0){80}}
\put(60,47){$\bullet$} %
\put(66,50){$b_{2}$}
%\put(35,47){$\bullet$}
\put(8,7){$\bullet$} %
\put(87,7){$\bullet$}%
\put(93,0){$a_{3}$}
\put(0,0){$a_{1}$}
\put(36,7){$\bullet$}
\put(33,0){$c_3$}
\put(63,7){$\bullet$}
\put(60,0){$b_{3}$}
\put(33,-20){(fig. 1)}
\end{picture}
\end{center}
\vskip 1cm
\vskip 3mm \noindent
{\bf Lemma 3.8}  $-r(b_{1}\vert a_{2},a_{3},b_{2},b_{3}) =
r'_{3}(a_{1},a_{2},a_{3},b_{1},b_{2},b_{3})$.
\vskip 3mm \noindent
{\bf Proof.}  The same as the one of lemma 3.8 in [G2]

Now let $\hat b_3$ be the  of the line $b_1b_2$ with the line
$a_1a_3$. Further, let $x$ be the intersection point of the lines $a_1b_2$
and $a_3b_1$. Let us denote by $c_3$
 the intersection point of the line $a_2 x$ with the line $a_1a_3$. Then
\begin{equation}
r'_{3}(a_{1},a_{2},a_{3},b_{1},b_{2},b_{3}) =  r(a_1,a_3,c_3,b_3)
\end{equation}
Indeed, by the well known theorem $r(a_1,a_3,\hat b_3,b_3) =-1$.

Now returning to a configuration $(v_1,...,v_6)$ (see fig 2)
\begin{center}
\begin{picture}(100,80)
\put(47,78){$a_{2}$}
\put(47,68){$\bullet$}
\put(50,70){\vector(2,-3){40}}
\put(50,70){\vector(-2,-3){40}}
\put(23,31){$\bullet$}
\put(12,33){$v_{1}$}
\put(71,31){$\bullet$}
\put(80,33){$v_{5}$}
\put(91,10){\vector(-1,0){80}}
\put(60,47){$\bullet$}
\put(66,50){$v_{2}$}
\put(25,50){$v_{4}$}
\put(35,47){$\bullet$}
\put(8,7){$\bullet$}
\put(87,7){$\bullet$}
\put(93,0){$a_{3}$}
\put(0,0){$a_{1}$}
\put(36,7){$\bullet$}
\put(33,0){$v_{6}$}
\put(63,7){$\bullet$}
\put(60,0){$v_{3}$}
\put(33,-20){(fig. 2)}
\end{picture}
\end{center}
\vskip 1cm
we see that one has proceed as follows: Put $b_1:=v_1,b_2:=v_2,b_3:=v_3$ and
apply the given above definition to the configuration
$(a_1,a_2,a_3,b_1,b_2,b_3)$
and then alternate. Notice that the configuration $(a_1,a_2,a_3,b_1,b_2,b_3)$
is defined by three flags $(v_1,v_1v_4),(v_2,v_2v_5),(v_3,v_3v_6)$.

\vskip 3mm \noindent
{\bf REFERENCES}
\begin{itemize}
\item[{[B1]}] Beilinson A.A.: {\it Higher regulators and values of
$L$-functions}, VINITI, 24 (1984), 181--238 (in Russian);
English translation: J. Soviet Math. 30 (1985), 2036--2070.
\item[{[B2]}] Beilinson A.A.: {\it Higher regulators for modular curves}
Contemporary Mathematics, vol. 55, 1987, 1-35.
\item[{[BL]}] Beilinson A.A., Levin A.M.: {\it Elliptic polylogarithms}.
Symposium in pure mathematics, 1994, vol 55, part 2, 101-156.
\item[{[Bl1]}] Bloch S.: {\it Higher regulators, algebraic $K$-
theory and zeta functions of elliptic curves}, Lect. Notes U.C.
Irvine, 1977.
\item[{[Bl2]}] Bloch S.: {\it 2 letters to Deninger
 regarding  [D1]} Fall 1990.
\item[{[Bl3]}] Bloch S.: {\it Lectures on algebraic cycles}, Duke Math.
Lect. Series, 1980.
\item[{[BMS]}] Beilinson A.A., MacPherson R.D. Schechtman V.V: {\it Notes
on motivic cohomology}. Duke Math. J., 1987 vol 55 p. 679-710
\item[{[GGL]}] Gabrielov A.M., Gelfand I.M., Losic M.V.:
{\it Combinatorial computation of characteristic classes}, Funct.\
Analysis and its Applications V. 9 No. 2 (1975) p. 103--115
and No. 3 (1975) p. 5--26 (in Russian).
\item[{[G1]}] Goncharov A.B.:{\it Geometry of configurations,
polylogarithms and motivic cohomology}.
Advances in Mathematics, 1995. 197 - 318.
\item[{[G2]}] Goncharov A.B., {\it Polylogarithms and motivic Galois
groups}, Symposium in pure mathematics, 1994, vol 55, part 2, p.  43 - 96.
\item[{[G3]}] Goncharov A.B., {\it Explicit construction of
characteristic classes}. Advances in Soviet mathematics, 1993,
vol 16, p. 169-210 (Special issue dedicated to I.M.Gelfand 80-th
birthday)
\item[{[G4]}] Goncharov A.B., {\it Special values of Hasse-Weil
$L$-functions and generalized Eisenstein-Kronecker
series}. To appear.
\item[{[GL]}] Goncharov A.B., Levin A.M. {\it Zagier's conjecture
    on $L(E,2)$}. Preprint IHES 1995.
\item[{[Del]}] Deligne P.: {\it Symbole modere} Publ. Math. IHES 1992
\item[{[D1]}] Deninger C.: {\it Higher order operations in Deligne
cohomology}. Inventiones Math. 122 N1 (1995).
\item[{[D2]}] Deninger C.: {\it Higher regulators and Hecke L-series of
imaginary quadratic fields I} Invent. Math. 96 (1989), 1-69.
\item[{[J]}] De Jeu R. {\it On $K_4^{(3)}$ of curves over number
    fields} Preprint 1995.
\item[{[S1]}] Suslin A.A.: {\it Homology of $GL_{n}$,
characteristic classes and Milnor's $K$-theory}.  Springer Lecture Notes
in Math. 1046 (1989), 357--375.
\item[{[S2]}] Suslin A.A.: {\it $K_{3}$ of a field and Bloch's
group}, Proceedings of the Steklov Institute of Mathematics 1991, Issue 4.
\item[{[Z1]}] Zagier D.:{\it Polylogarithms, Dedekind zeta
functions and the algebraic $K$-theory of fields}. Arithmetic
Algebraic Geometry (G.v.d.Geer, F.Oort, J.Steenbrink, eds.),
Prog. Math., Vol 89, Birkhauser, Boston, 1991, pp. 391--430.
\item[{[Z2]}] Zagier D.:{\it The Bloch - Wigner -Ramakrishnan
polylogarithm function} Math. Ann. 286 (1990), 613-624

\end{itemize}

\end{document}